\documentclass[sn-basic]{sn-jnl}
\jyear{2023}%

\theoremstyle{thmstyleone}%
\theoremstyle{thmstyletwo}%
\theoremstyle{thmstylethree}%

\allowdisplaybreaks

\usepackage{mathtools,enumitem}
\usepackage[capitalise,nameinlink,noabbrev]{cleveref}
\usepackage[nolist,printonlyused,withpage]{acronym}
\usepackage[binary-units=true,exponent-product = \cdot, output-product = \cdot]{siunitx}
\usepackage{bm}

\numberwithin{equation}{section}
\crefname{subsection}{Subsection}{Subsections}
\crefformat{equation}{(#2#1#3)}
\crefformat{align}{(#2#1#3)}
\crefformat{enumi}{(#2#1#3)}
\crefname{figure}{Figure}{Figures}

\newcommand{\eps}{\varepsilon}

\newcommand{\ECM}{{E\hspace{-.05em}C\hspace{-.15em}M}}
\newcommand{\MDE}{{M\hspace{-.15em}D\hspace{-.1em}E}}
\newcommand{\TAF}{{T\hspace{-.15em}A\hspace{-.1em}F}}
\newcommand{\CMT}{{C\hspace{-.15em}M\hspace{-.1em}T}}
\newcommand{\ecm}{\phi_{\ECM}}
\newcommand{\cmt}{\phi_{\CMT}}
\newcommand{\taf}{\phi_{\TAF}}
\newcommand{\mde}{\phi_{\MDE}}

\renewcommand{\rel}{\textup{rel}}
\newcommand{\dd}{\textup{d}}
\newcommand{\dx}{\,\textup{d}x}

\newcommand{\dt}{\,\textup{d}t}
\newcommand{\ds}{\,\textup{d}s}

\renewcommand{\A}{\mathbb{A}}
\renewcommand{\R}{\mathbb{R}}
\newcommand{\N}{\mathbb{N}}
\renewcommand{\p}{\partial}
\newcommand{\pt}{\p_t}

\newcommand{\ov}{\overline}

\newcommand{\apo}{\textup{apo}}
\newcommand{\pro}{\textup{pro}}
\renewcommand{\deg}{\textup{deg}}
\renewcommand{\kill}{\textup{kill}}
\renewcommand{\div}{\textup{div}}
\newcommand{\tr}{\textup{tr}}

\newcommand{\calE}{\mathcal{E}}

\newcommand{\calH}{\mathcal{H}}

\newcommand{\scrL}{\mathscr{L}}

\newcommand{\CH}{{\mathbb{C}\mathbb{H}}}
\newcommand{\RD}{{\mathbb{R}\mathbb{D}}}
\renewcommand{\OD}{{\mathbb{O}\mathbb{D}}}

\DeclareFontFamily{U}{bbold}{}
\DeclareFontShape{U}{bbold}{m}{n}
{
	<-5.5> s*[1.069] bbold5
	<5.5-6.5> s*[1.069] bbold6
	<6.5-7.5> s*[1.069] bbold7
	<7.5-8.5> s*[1.069] bbold8
	<8.5-9.5> s*[1.069] bbold9
	<9.5-11> s*[1.069] bbold10
	<11-15> s*[1.069] bbold12
	<15-> s*[1.069] bbold17
}{}

\DeclareRobustCommand{\bbI}{%
	\text{\usefont{U}{bbold}{m}{n}1}%
}

\def\signed #1{{\leavevmode\unskip\nobreak\hfil\penalty50\hskip2em
		\hbox{}\nobreak\hfil(#1)%
		\parfillskip=0pt \finalhyphendemerits=0 \endgraf}}

\newsavebox\mybox

\acrodef{cmt}[CMT]{chemotherapy}
\acrodef{ctrw}[CTRW]{continuous time random walk}
\acrodef{dfb}[DFB]{Darcy--Forchheimer--Brinkman}
\acrodef{ecm}[ECM]{extracellular matrix}
\acrodef{fem}[FEM]{finite element method}
\acrodef{fde}[FDE]{fractional differential equation}
\acrodef{hpc}[HPC]{high performance computing}
\acrodef{mde}[MDE]{matrix-degrading enzyme}
\acrodef{msd}[MSD]{mean square displacement}
\acrodef{ode}[ODE]{ordinary differential equation}
\acrodef{pde}[PDE]{partial differential equation}
\acrodef{pdf}[PDF]{probability density function}
\acrodef{rde}[RDE]{reaction-diffusion equation}
\acrodef{taf}[TAF]{tumor angiogenesis factor}
\acrodef{vegf}[VEGF]{vascular endothelial growth factor}
\acrodef{vgm}[VGM]{vascular graph method}

\begin{document}
\title[Tumor evolution models of phase-field type with nonlocal effects]{\mbox{}
\vspace{-3cm} \\ Tumor evolution models of phase-field type with nonlocal effects and angiogenesis}
    \author{\fnm{Marvin} \sur{Fritz}}\email{marvin.fritz@ricam.oeaw.ac.at}
    \affil{\orgdiv{Computational Methods for PDEs}, \orgname{Johann Radon Institute for Computational and Applied Mathematics}, \orgaddress{\city{Linz}, \country{Austria}}}
	
	\abstract{
		In this survey article, a variety of systems modeling tumor growth  are discussed. In accordance with the hallmarks of cancer, the described models incorporate the primary characteristics of cancer evolution. Specifically, we focus on diffusive interface models and follow the phase-field approach that describes the tumor as a collection of cells. Such systems are based on a multiphase approach that employs constitutive laws and balance laws for individual constituents. In mathematical oncology, numerous biological phenomena are involved, including temporal and spatial nonlocal effects, complex nonlinearities, stochasticity, and mixed-dimensional couplings. Using the models, for instance, we can express angiogenesis and cell-to-matrix adhesion effects. Finally, we offer some methods for numerically approximating the models and show simulations of the tumor's evolution in response to various biological effects.}


    \maketitle
	
	\section{Introduction}

Cancer is among the main global causes of death.
According to \cite{sung2021global}, there were 19.3 million new cancer diagnoses and 9.96 million cancer-related deaths worldwide.
By 2040, the yearly number of new cancer cases is projected to reach 30.2 million, with 16.3 million fatalities attributable to cancer.
Each tumor is distinct and dependent on a variety of characteristics.
There is no guaranteed procedure for curing cancer, nor is its cause entirely known.
Utilizing mathematical models to precisely depict tumor progression is the primary objective of mathematical oncology. 

The key hallmarks of cancer evolution are described by \cite{hanahan2000hallmarks,hanahan2011hallmarks}
and for mathematical oncology to be successful, these characteristics should be met. As a primary advantage of a realistic mathematical model, cancer progression can be forecasted and physicians will be able to simply press a button on their computers to initiate a simulation portraying the patient's tumor and its development.
Ideally, this process is combined with a focused therapy that improves the cancer's prognosis. However, one must first guarantee that the model is well-posed, both mathematically and in terms of accurately representing the movement of actual cancer.
The second point can only be investigated using data and model verification through prediction; see the survey article by \cite{oden_2018} for more information on this topic.
The direction of this survey paper is toward the first point.
We must ensure that these models are mathematically valid, have a solution, and that nothing nonsensical occurs.
Then, one can consider a numerical strategy for the model that will provide a rapid, accurate, and stable representation of the tumor's evolution on the physician's monitor. 

There is an abundance of literature on the mathematical modeling of tumor evolution, which is a positive development.
Different groups develop distinct models and procedures and with this diversification, it is hoped that researchers will be able to accurately forecast the progression of malignancies.
In describing the phenomena of the world, partial differential equations (PDEs) are ubiquitous; they model the flow of liquids and gases (Navier--Stokes equations), the evolution of a quantum state (Schr\"odinger equation), thermal conduction (heat equation), spinodal decomposition (Cahn--Hilliard equation), and many others.
Complicated models may include nonlinearities, temporal and spatial nonlocalities, and mixed-dimensional couplings in response to complex processes.

Initially, tumor models were expressed as a free boundary problem.
We refer to \cite{greenspan1976growth}, which treated the tissue as a porous media and calculated the convective velocity field using Darcy's law.
Such models have been expanded upon in various works, and we direct you to the previous reviews by \cite{bellomo2000modeling} and \cite{roose2007mathematical}.
Since then, numerous distinct models have been formulated and in particular, we follow the path of diffusive interface models in which the tumor is characterized as a collection of cells using a fourth-order PDE.
These models are based on a multiphase method employing constitutive laws, thermodynamic principles, and balance rules for single constituents, which dates back to the works of \cite{cristini2003nonlinear}, \cite{cristini2010multiscale}, \cite{frieboes2010three}, and \cite{wise2008three} starting in 2003.

This work is organized as follows:  In \cref{Sec:Modeling}, we examine tumor evolution models and follow a technique based on continuum mixture theory. In this regard, we present the Cahn--Hilliard equation, the fundamental model of our tumor growth systems.
We provide a multiphase tumor growth model consisting of numerous components and biological processes.
In particular, we investigate the effects of the extracellular matrix, tumor cell stratification, the release of matrix degenerating enzymes and tumor angiogenesis factors, stochasticity, mechanical deformation, chemotherapeutic influence, memory effects, subdiffusion, and nonlocal phenomena including cell-to-cell adhesion and cell-to-matrix adhesion.  Further, we highlight each phenomena by numerical simulations and illustrations. We state the ideas for the numerical approximations of the introduced models in \cref{Sec:Numerics}.

    \section{Modeling of Tumor Growth}\label{Sec:Modeling}

\noindent We propose mathematical oncology models that abstract a number of the known significant mechanisms involved in tumor growth, decline, and therapeutic therapy in real tissue.
The systems are designed to reflect mesoscale and macroscale dynamics, with fields representing volume fractions of mass concentrations of diverse species that determine tumor composition.
Several authors, including \cite{araujo2004history},  \cite{fritz2022well}, \cite{garcke2016cahn,garcke2018multiphase}, \cite{lima2014hybrid} and \cite{wise2008three}, have produced localized versions of multiphase models over the past decade. Balance laws of continuum mixture theory are used to derive the model equation, see also \cite{byrne2003modelling}, \cite{cristini2009nonlinear}, and \cite{oden2016toward, oden2010general}.

In \cref{Sec:Math:Proto}, we present the prototype system for modeling tumor growth -- the Cahn--Hilliard equation with concentration-dependent mobility. In a generic framework, we provide in \cref{Sec:Mod:Mult} a multiple constituent model derived from the mass balance law and a Ginzburg--Landau type energy.
As an illustration, we provide the four-species model developed by \cite{hawkins2012numerical}.
In \cref{Sec:Mod:ECM}, we incorporate stratification and invasion due to ECM deterioration into the model. In the following subsections, additional biological phenomena will be added to the stratified tumor model.  We incorporate spatial and temporal nonlocalities in \cref{Sec:Mod:Nonlocal}, stochasticity by a cylindrical Wiener process in \cref{Sec:Mod:Stochastic}, mechanical deformation in \cref{Sec:Mod:Mechanical}, chemotherapeutic influence in \cref{Sec:Mod:Chemo}, and lastly, angiogenesis in mixed-dimensional couplings in \cref{Mod:Mixed}.

\subsection{Prototype model: The Cahn--Hilliard equation} \label{Sec:Math:Proto}
The Cahn--Hilliard equation is the prototypical model for tumor growth.
It is a phase-field equation of the diffuse-interface type, and it possesses the essential attribute of having a solution that is either 0 or 1, or a smooth transition phase in between.
We define the 1-phase as the manifestation of tumor cells, whereas the 0-phase represents the absence of malignant cells. 

Let $\phi_1$ and $\phi_2$ represent the concentrations of two components, and it holds $\phi_1+\phi_2=1$. This indicates that the concentrations describe local portions, such as those found in binary alloys.
They comply with the mass conservation law 
	\begin{equation*}
	\p_t \phi_i = - \div J_i, \quad i \in \{1,2\},
	\end{equation*}
	where the mass flux of the $i$-th component is denotes by $J_i$. We assume that the fluxes fulfill the condition $J_1+J_2=0$ and we reduce the equations by defining the quantities $\phi=\phi_1-\phi_2$ and $J=J_1-J_2$, which yields
	$$  \p_t \phi = - \div J.$$
	Here, the flux $J$ is given by the negative of the gradient of the chemical potential $\mu$, i.e., $J=-\nabla \mu$. In \cite{gurtin1996generalized}, a mechanical version of the second law of thermodynamics was introduced
by providing an augmented mass flux $J=-m(\phi)\nabla\mu$ with some mobility function $m$ for describing microscopic interactions.
	Following \cite{cahn1958free}, the chemical potential $\mu$ is given by the G\^ateaux derivative of the Ginzburg--Landau energy functional
	\begin{equation} \label{Eq:Ginzburg}
	\calE(\phi) = \int_\Omega \Big\{ \Psi(\phi) + \frac{\eps^2}{2} \vert \nabla \phi\vert ^2 \Big\} \, \dd x.\end{equation}
	Here, the parameter $\eps$ expresses the interfacial width and $\Psi$ describes a double-well potential with zeros at $0$ and $1$, e.g., the Landau potential
	$
	\Psi(\phi)=\frac14 \phi^2 (1-\phi)^2. 
	$
 Hence, the Cahn--Hilliard equation with concentration-dependent mobility reads
\begin{equation} \label{Eq:Cahn} \boxed{
		\begin{gathered}
		\textbf{Cahn--Hilliard equation}  \\
		\begin{aligned}
	\p_t \phi ={}& \div(m(\phi) \nabla \mu) \\  \mu ={}& \Psi'(\phi) - \eps^2 \Delta \phi \end{aligned}\end{gathered}}
\end{equation}
	  Usually, the mobility function takes the form $m(\phi)=M\phi^2(1-\phi)^2$ for some $M>0$. The scenario of constant mobility has been exhaustively examined, and well-posedness can be demonstrated through the use of sufficient assumptions, as done in \cite{miranville2019cahn}. A proof or counterexample of uniqueness in the case of degenerate mobility remains unsolved for the class of degenerate fourth-order parabolic equations.

\subsection{Base system: Multiple constituent model} \label{Sec:Mod:Mult}
Multiple mechanical and chemical species can coexist at a given place $x$ in a given domain $\Omega \subset \R^d$, $d \in \N$, within the continuum mixture theory paradigm.
For a medium with $N$ interacting constituents, the volume fraction of each species is therefore represented by a field $\phi_\alpha$, $1\leq \alpha \leq N$, with value $\phi_\alpha(t,x)$ at $x\in \Omega$, and time $t\geq 0$.
For convenience, we compile the model's components in the following $N$-tuple
$$
 \phi_\A =  (\phi_\alpha)_{\alpha \in \A},
$$
where $\A$ is an index set that is further disjointly separated between the phase-field index set $\CH$, the reaction-diffusion indices $\RD$, and the evolution indices $\OD$ that correspond to abstract ordinary differential equations (ODEs).

Following \cite{lima2014hybrid,lima2015analysis}, the constituents $\phi_\alpha$, $\alpha \in \A$, are governed by the  extended mass balance law
\begin{equation} \label{Eq:MassBalance}
\p_t \phi_\alpha+\text{div}(\phi_\alpha v_\alpha)=- \text{div} J_\alpha(\phi_\A) +S_\alpha(\phi_\A).
\end{equation}
Here, $v_\alpha$ is the cell velocity of the $\alpha$-th  constituent, and $S_\alpha$ is a species-dependent mass source term.
We refer to the system as closed if it holds $\sum_{\alpha \in \A} S_\alpha(\phi_\A)=0$.
In addition, $J_\alpha$ represents the flux of the $\alpha$-th constituent, which is proportional to the negative gradient of the chemical potential multiplied by a mobility function
\begin{equation} \label{Eq:Flux}
J_\alpha(\phi_\A) = - m_\alpha(\phi_\A) \nabla \mu_\alpha.
\end{equation}
Here, $\mu_\alpha$ represents the chemical potential of the $\alpha$-th species, and $m_\alpha$ represents the mobility function, which may depend on all constituents.
In our applications, we typically take the mobilities
\begin{equation} \label{Eq:Mobility} \begin{aligned} m_\alpha(\phi_\A) ={}&M_\alpha \phi_\alpha^2 (1-\phi_\alpha)^2, && \alpha \in \CH, \\ m_\beta(\phi_\A) ={}& M_\beta, && \beta \in \RD, \\
m_{\gamma}(\phi_\A) ={}&0, && \gamma \in \OD,&&
\end{aligned}\end{equation} where $M_\alpha>0$ are constants.  Similarly to the prototype model, see \cref{Sec:Math:Proto}, we define the chemical potential $\mu_\alpha$ as the G\^ateaux derivative of the Ginzburg--Landau energy with respect to $\phi_\alpha$. We propose the system's energy
\begin{equation}\mathcal{E}(\phi_\A)= \int_\Omega \Big\{ \Psi(\phi_{\CH}) + \Phi(\phi_\A)+ \sum_{\alpha \in \CH}  \frac{\eps_\alpha^2}{2} \lvert\nabla \phi_\alpha\rvert^2  + \sum_{\beta \in \RD} \frac{D_\beta}{2} \phi_\beta^2 \Big\}  \text{ d}x,
\label{Eq:GinzburgLandau}
\end{equation}
where $\varepsilon_\alpha$, $\alpha \in \CH$, is a parameter related to the thickness of the contact separating the various cell kinds. As we will see later, the function $\Phi$ explains adhesion mechanisms such as chemotaxis and haptotaxis.
Finally, $\Psi$ represents a double-well potential as in the generic Cahn--Hilliard equation \cref{Eq:Cahn}, e.g., it may be of Landau type, for which we list two alternatives
$$\begin{aligned} \Psi(\phi_\CH)&={}C_{\Psi} \bigg(\sum_{\alpha \in \CH} \phi_\alpha\bigg)^2 \bigg(1-\sum_{\alpha \in \CH} \phi_\alpha \bigg)^2, \\
\Psi(\phi_\CH) &={} \sum_{\alpha \in \CH} C_{\Psi_\alpha} \phi_\alpha^2(1-\phi_\alpha)^2, \end{aligned}
$$
where $C_\Psi$ and  $C_{\Psi_\alpha}$ are given prefactors. As another possibility, we could select a logarithmic potential of Flory--Huggins type, see \cite{cherfils2011cahn} and \cite{frigeri2018on}.

We calculate the G\^{a}teaux derivatives of the Ginzburg--Landau energy \cref{Eq:GinzburgLandau} with respect to the stated constituents and therefore, the corresponding chemical potentials read
\begin{equation*} %
\begin{aligned}
\mu_\alpha ={}& \p_{\phi_\alpha} \Psi(\phi_\CH) + \p_{\phi_\alpha} \Phi(\phi_\A) - \eps^2_\alpha \Delta \phi_\alpha, && \alpha \in \CH , \\
\mu_\beta ={}& D_\beta \phi_\beta  + \p_{\phi_\beta} \Phi(\phi_\A), && \beta \in \RD , \\
\mu_\gamma ={}& \p_{\phi_\gamma} \Phi(\phi_\A), &&\gamma \in \OD,
\end{aligned}
\end{equation*}
and combining the chemical potentials with the mass balance laws \cref{Eq:MassBalance}--\cref{Eq:Mobility}, it yields the multispecies model:
\begin{equation} \label{Eq:MultipleGeneral} \boxed{ 
	\begin{gathered} 
	\textbf{Multiple constituent model} \\
	\begin{aligned}
\pt \phi_\alpha\!+\! \div(\phi_\alpha v_\alpha) ={}& \div \big(M_\alpha\phi_\alpha^2(1-\phi_\alpha)^2 \nabla \mu_\alpha\big) + S_\alpha( \phi_\A) &&\alpha \in \CH   \\
\mu_\alpha ={}& \p_{\phi_\alpha} \Psi(\phi_\CH)+ \p_{\phi_\alpha} \Phi(\phi_\A) - \eps^2_\alpha \Delta \phi_\alpha  &&\alpha \in \CH \\
\pt \phi_\beta \!+\!    \div(\phi_\beta v_\beta)                     ={}& \div\big(M_\beta \nabla \big( D_{\beta} \phi_\beta \!+\! \p_{\phi_\beta} \Phi(\phi_\A) \big) \big)\! +\!S_{\beta}( \phi_\A) &&\beta \in \RD \\
\pt \phi_\gamma                        ={}& S_{\gamma}( \phi_\A) &&\gamma \in \OD                     
\end{aligned} \end{gathered}} \end{equation}

\subsubsection{Four-species tumor growth model}
We begin with a straightforward illustration of a tumor growth model based on the suggested multiple constituent system \cref{Eq:MultipleGeneral}.

The article by \cite{hawkins2012numerical} presents the most fundamental model of tumor growth, which forms the basis of this theory.
The volume fractions of cancer cells, healthy cells, nutrient-rich extracellular water, and nutrient-poor extracellular water were taken into account.
Such a system is referred to as the “four-species model,” and \cite{garcke2016global,garcke2017analysis,garcke2017well} investigated the model's mathematical well-posedness.
In addition, we cite \cite{frigeri2015diffuse,frigeri2017diffuse} for an examination of degenerating mobility functions.
Due to the fact that the model is based on a fourth-order PDE with concentration-dependent mobilities, even for the prototype model \cref{Eq:Cahn}, the uniqueness of weak solutions is unresolved; see \cite{elliott1996cahn} for more information.
\cite{colli2017optimal} investigated the four-species model in relation to an optimal control problem, whereas \cite{miranville2019long} and \cite{cavaterra2011cahn} investigated the long-term behavior of the solution.
Various velocity models have been introduced to the four-species model to account for fluid movement in the progression of cancer.
The cells are represented as viscous, inertia-free fluids, and the fluid mixture's velocity is modeled in a volume-averaged sense.
Such an assumption is reasonable, given that the cells are densely packed.
\cite{garcke2016cahn} modeled the velocity by the Darcy law in the four-species model, and \cite{garcke2018cahn} examined this model analytically.
In \cite{ebenbeck2019analysis,ebenbeck2019cahn} and in \cite{fritz2019unsteady}, this law was extended to the Darcy--Brinkman equation and the time-dependent Darcy--Forchheimer--Brinkman equation, respectively.
Authors have also approximated the velocity as a Stokes flow (see \cite{franks2003interactions} and \cite{friedman2006free,friedman2016free}), and the Darcy--Brinkman equation can be viewed as an interpolation between Darcy and Stokes flow.
The inclusion of a velocity equation in a Cahn--Hilliard system is not innovative in and of itself, as it has been done by \cite{lee2002modeling} without the application to tumor growth.
These strategies have been modified to accommodate the new system, which incorporates nontrivial effects such as chemotaxis, proliferation, and nonlinear source functions. 

\begin{figure}[htb!]
    \centering
    \includegraphics[width=.8\textwidth]{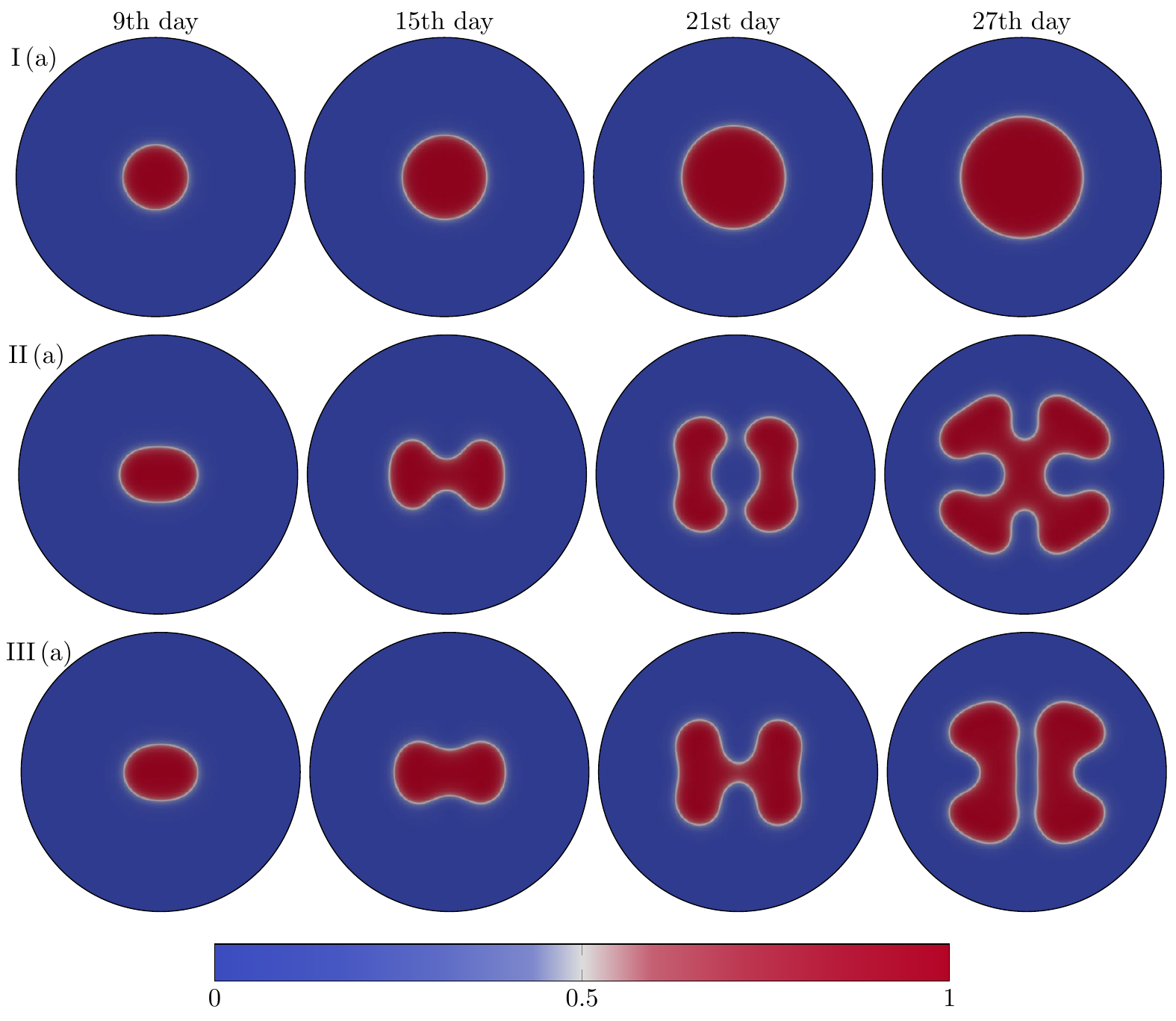}
    \caption{Evolution of tumor mass $\phi_T$ with a slightly elliptic initial condition on the 9th, 15th, 21st and 27th day; we present three different variations of the model: I. without velocity, II. unsteady Darcy--Brinkman law, III. unsteady Darcy--Forchheimer--Brinkman law; figure taken with permission from Figure 7 in \cite{fritz2019unsteady}.}
    \label{Fig:Image_Flow}
\end{figure}

We choose $\vert\A\vert=2$ constituents and set $\A=\{T,\sigma\}$, $\CH=\{T\}$, $\RD=\{\sigma\}$, and $\OD=\emptyset$. It is understood that the volume fraction of tumor cells $\phi_T$ represents an averaged cell concentration, a homogenized representation of several thousands of cells.
Field $\phi_\sigma$ is representative of the local nutrient content.
In addition, we present the adhesion function $\Phi(\phi_T,\phi_\sigma)=-\chi_c\phi_T \phi_\sigma$ in energy \cref{Eq:GinzburgLandau} for a particular chemotaxis parameter $\chi_c>0$.
For the tumor cells and the nutrients, we assume a volume-averaged velocity.
This assumption of a volume-averaged velocity is fair given the dense packing of the cells.
When all the assumptions are inserted into the multispecies model, the result is the so-called four-species model. 
\begin{equation} \label{Eq:FourSpecies} \boxed{ 
	\begin{gathered} 
	\textbf{Four-species model} \\
	\begin{aligned}
		\pt \phi_T+ \div(\phi_T v) ={}& \div \big(M_T\phi_T^2(1-\phi_T)^2  \nabla \mu_T\big) + S_T( \phi_T,\phi_\sigma)    \\
		\mu_T ={}& \Psi'(\phi_T) - \chi_c \phi_\sigma- \eps_T^2 \Delta \phi_T \\
		\pt \phi_\sigma +  \div(\phi_\sigma v)                     ={}& \div\big(M_\sigma \nabla( D_{\sigma} \phi_\sigma -\chi_c\phi_T  ) \big) +S_{\sigma}( \phi_T,\phi_\sigma) 
\end{aligned} \end{gathered} } \end{equation}
In the case of an absent velocity $v=0$, this model is studied in \cite{garcke2017analysis,garcke2017well} with respect to the existence of weak solutions. If the flow is governed by Darcy's law $$\begin{aligned} v&=-K\nabla p + S_v(\phi_T,\phi_\sigma), \\ \div \, v &= 0,\end{aligned}$$ then we refer to \cite{garcke2016cahn} and \cite{garcke2016global}. The pressure is denoted by $p$, the permeability factor by $K>0$, and $S_v$ is called the Korteweg force \cite{frigeri2018on}. Alternatively, the flow has been governed by the Brinkman law \citep{ebenbeck2019analysis,ebenbeck2019cahn}, the unsteady Darcy--Forchheimer--Brinkman law \citep{fritz2019unsteady}, and the Navier--Stokes equations \citep{lam2017thermodynamically,he2021global} in literature. Numerically, we present a comparison of different flow models and their influence in the four-species model, see \cref{Fig:Image_Flow}. We refer to \cref{Sec:Numerics} below for further details on the techniques for discretizing the PDEs in time and space. We notice that the flow is highly influential on the evolution of the tumor by drastically changing the growth directions of the tumor mass.

Source functions that are expressed as sink and source terms are of particular importance.
Tumors absorb the nutrients; hence, tumor growth is proportional to nutrient depletion.
In addition, programmed cell death (also known as apoptosis) occurs, and these dead cells become nutrients.
Consequently, we consider the source function
$$S_T(\phi_T,\phi_\sigma)=-S_\sigma(\phi_T,\phi_\sigma)=\lambda^\pro_T \phi_\sigma\phi_T(1-\phi_T) - \lambda^\apo_T \phi_T,$$
where $\lambda_T^\pro$ is called proliferation rate and $\lambda_T^\apo$ apoptosis rate.

The system \cref{Eq:FourSpecies} is also referred to as the “four-species model,” \citep{hawkins2012numerical,oden2010general,lima2014hybrid} because it can be derived from four constituents: the volume fraction of tumor cells $\phi_T$, healthy cells $\phi_C$, nutrient-rich extracellular water $\phi_\sigma$, and its nutrient-poor counterpart $\phi_{\sigma_0}$.
Consequently, the four variables are governed by the law of mass balance, see \cref{Eq:MassBalance}, for $\A=\{T,C,\sigma,\sigma_0\}$. One sets $\phi_T=1-\phi_C$ and $\phi_\sigma=1-\phi_{\sigma_0}$. Thus, one can eliminate the superfluous constituents $\phi_C$ and $\phi_{\sigma_0}$ from the system and obtains the four-species system \cref{Eq:FourSpecies}.

\subsection{Phase separation in an ECM} \label{Sec:Mod:ECM}
The “microenvironment” of a solid tumor is a patch of vascularized tissue in a living subject, such as within an organ, that contains a colony of tumor cells and other components.
The tumor is contained within an open-bounded region $\Omega \subset \R^3$ and is supported by a network of collagen, enzymes, and other proteins that comprise the extracellular matrix (ECM).
We are focusing on developing phenomenological descriptions of tumor cell colony growth that capture both mesoscale and macroscale phenomena.

When tumor cells endure hypoxia or necrosis, these four-species models are inadequate for representing the formation of an early tumor whose evolution is primarily determined by proliferation.
Indeed, a larger and more advanced tumor tends to become stratified \citep{roose2007mathematical}, meaning that the tumor tissue is subdivided into numerous layers, each with its own properties.
Typically, tumors are separated into three phases: \medskip
\begin{itemize} \itemsep0.5em
	\item Rapidly proliferating outer rim.
	\item Intermediate quiescent layer with cells suffering from hypoxia.
	\item Necrotic core with perished cells. \medskip
\end{itemize}
	Multiphase models with multiple cell species and nutrients have been studied in the works \cite{wise2008three}, \cite{escher2011analysis}, \cite{sciume2014three}, \cite{garcke2018multiphase}, \cite{araujo2004history}, \cite{astanin2008multiphase}, \cite{frieboes2010three}, \cite{frigeri2018on}, \cite{dai2017analysis}, and \cite{fritz2019local,fritz2021analysis,fritz2021modeling}.
In the hypoxic phase, tumor cells are quiescent and release matrix-degrading enzymes (MDEs), which degrade the ECM and allow nutrients to flow.
This procedure allows tumor cells to move into the tissue and is the initial stage in simulating metastasis.
Simply put, the ECM works as a wall that regulates the flow of nutrients around the tumor.
Several authors \citep{chaplain2011mathematical,engwer2017structured,stinner2014global,sfakianakis2020hybrid,shuttleworth2020cell,sciume2014tumor} have examined the ECM in reaction-diffusion type tumor models.
We investigated the ECM in a Cahn--Hilliard type model \citep{fritz2019local}, and it was also included in our successive research \citep{fritz2021analysis,fritz2021modeling}.

The field of the tumor cells $\phi_T$ can be represented by the sum $$\phi_T=\phi_P+\phi_H+\phi_N,$$ of the three components $\phi_P$, $\phi_H$, $\phi_N$ that describe the volume fractions of the proliferative, hypoxic, and necrotic cells, respectively. They are characterized by: \medskip

 \begin{itemize} \itemsep.5em
 	\item Proliferative cells $\phi_P$ are those with a high probability of undergoing mitosis, dividing into twin cells, and fostering tumor growth. 
 	\item Hypoxic cells $\phi_H$ are tumor cells that lack sufficient resources, such as oxygen, to proliferate or continue to proliferate. 
 	\item Necrotic cells $\phi_N$ have died owing to nutrient deficiency. \medskip
 	\end{itemize} 
  
\noindent In response to hypoxia, tumor cells produce an enzyme that promotes cell motility and stimulates the secretion of angiogenesis-stimulating substances $\phi_\TAF$.
The most frequently mentioned of these substances is vascular endothelial growth factor (VEGF), which induces endothelial cells to proliferate and form the tubular shape of blood vessels, which then extend to form new arteries that supply nutrition to hypoxic cells.

\begin{figure}[htb!]
    \centering
    \includegraphics[width=.245\textwidth]{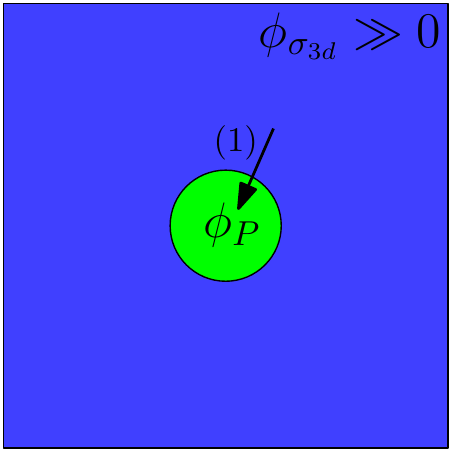}\!
\includegraphics[width=.245\textwidth]{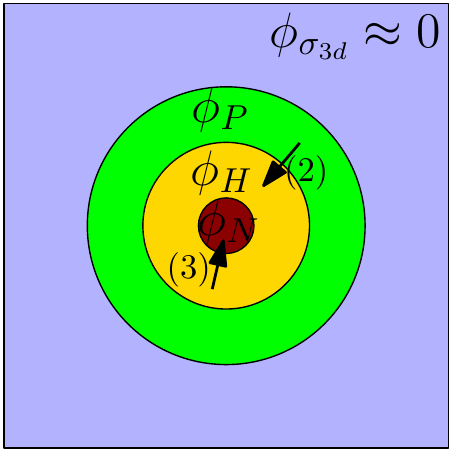}\!
\includegraphics[width=.245\textwidth]{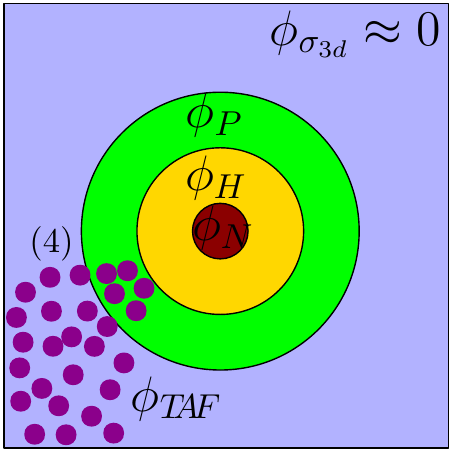}\!
\includegraphics[width=.245\textwidth]{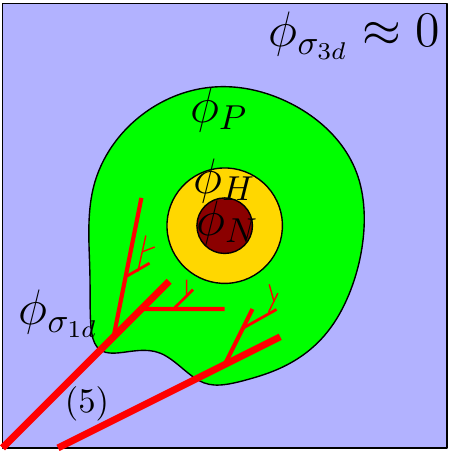}
    \caption{Depiction of angiogenesis and growth of capillaries after the proliferative tumor phase becomes hypoxic due to nutrient shortage.}
    \label{Fig:Separation}
\end{figure}

In addition, hypoxic cells release MDEs such as urokinase-plasminogen and matrix metalloproteinases, as indicated by the volume fraction $\phi_\MDE$, which erode the ECM, whose density is represented by $\phi_\ECM$.
This procedure permits tumor cells $\phi_T$ to infiltrate, hence increasing the number of tumor cells in the ECM domain and the probability of metastasis.
The following is a simplified explanation of the impacts of the tumor's evolution and it is also depicted in \cref{Fig:Separation}. \medskip

\begin{enumerate} \itemsep.5em
	\item[(1)] Outer proliferative layer absorbs nutrients and expands ($\phi_P\!\!\uparrow$,  $\phi_{\sigma}\!\downarrow$).
	\item[(2)] Inner tumor layer changes to hypoxic ($\phi_H\!\uparrow$).
	\item[(3)] Tumor core changes to necrotic ($\phi_N \uparrow$).
	\item[(4)] Hypoxic cells send out MDEs and TAFs ($\taf\!\uparrow$, $\mde\!\uparrow$).
	\item[(5)] TAFs trigger angiogenesis and initiate the sprouting of vessels  ($\phi_H\!\!\downarrow, \phi_P\!\uparrow$), \\[0.1cm] and MDEs erode the ECM, i.e., tumor cells migrate  ($\ecm\!\downarrow$, $\phi_H\!\downarrow, \phi_P\!\uparrow$). \medskip
\end{enumerate}

\noindent We collect the constituents within the following tuple:
$$
\phi_\A = (\phi_\alpha)_{\alpha \in \A} = ( \phi_P, \phi_H, \phi_N, \phi_\sigma, \ecm, \mde, \taf  ),
$$
with $\A=\{P,H,N,\sigma,\ECM,\MDE,\TAF\}$. We differentiate between the tumor phase-field indices $\CH=\{P,H,N\}$, the reaction-diffusion indices $\RD=\{\sigma,\MDE,\TAF\}$, and the evolution index set $\OD=\{\ECM\}$ using the setup of the multiple constituent model \cref{Eq:MultipleGeneral} in \cref{Sec:Mod:Mult}.
The necrotic cells are immobile and only gain mass from the hypoxic cells, which lack nutrients.
Therefore, the necrotic cells' mobility is set to zero, i.e., it holds $m_N=v_N=0$.
Still, necrotic cells are counted as a phase-field variable and constitute a component of $\CH$ rather than the ODEs because they influence the double-well potential and inherit their phase-field structure from the hypoxic phase-field variable.
Assuming that haptotaxis and chemotaxis are part of the system, we calculate the adhesion force $$\Phi(\phi_\A)=-(\phi_P+\phi_H)(\chi_c \phi_\sigma + \chi_h \phi_\ECM),$$ 
where $\chi_c$ and $\chi_h$ are the chemotaxis and haptotaxis components, respectively.
The adhesion force only operates on live (proliferative and hypoxic) cells, while necrotic cells are excluded from this process.
Consequently, the equations for the phase-field variables $(\phi_\alpha)_{\alpha \in \CH}$ are derived from the multiple constituent model \cref{Eq:MultipleGeneral} and read as follows:
\begin{equation} \label{Eq:Stratified} \boxed{
	\begin{gathered} 
	\textbf{Stratified tumor growth model with ECM: $\CH$} \\
	\begin{aligned}
		\pt \phi_P+ \div(\phi_P v) ={}& \div \big(M_P \phi_P^2(1-\phi_P)^2 \nabla \mu_P\big) + S_P( \phi_\A)    \\
		\mu_P                      ={}&   \p_{\phi_P} \Psi(\phi_\CH) - \eps^2_P \Delta \phi_P - \chi_c \phi_\sigma-\chi_h \ecm \\
		\pt \phi_H+ \div(\phi_H v) ={}& \div \big(M_H \phi_H^2(1-\phi_H)^2 \nabla \mu_H\big)+S_H( \phi_\A)      \\
		\mu_H                      ={}&   \p_{\phi_H} \Psi(\phi_\CH) - \eps^2_H \Delta \phi_H - \chi_c \phi_\sigma-\chi_h \ecm \\
		\pt \phi_N                 ={}& S_N( \phi_\A)
	\end{aligned} \end{gathered} }\end{equation}
We assume a volume-averaged velocity $v=v_\alpha$ for the fields $\phi_P$, $\phi_H$ and $\phi_\sigma$, which shall be governed by Darcy law for the sake of simplicity. Moreover, we consider the following source functions
$$
	\begin{aligned}
S_P(\phi_\A) ={}& \lambda_P^\pro \phi_\sigma \phi_P(1- \phi_T) - \lambda_P^\apo \phi_P - \lambda_{P\!H}  \calH(\sigma_{P\!H} - \phi_\sigma)\phi_P \\& + \lambda_{H\!P} \calH(\phi_\sigma - \sigma_{H\!P}) \phi_H,  \\
S_H(\phi_\A) ={}&  \lambda_{H}^\pro \phi_\sigma \phi_H (1-\phi_T)-\lambda_H^\apo \phi_H + \lambda_{P\!H}  \calH(\sigma_{P\!H} - \phi_\sigma) \phi_P  \\& - \lambda_{H\!P} \calH(\phi_\sigma - \sigma_{H\!P}) \phi_H - \lambda_{H\!N} \calH(\sigma_{H\!N} - \phi_\sigma)\phi_H,  \\
S_N(\phi_\A) ={}& \lambda_{H\!N} \calH(\sigma_{H\!N} - \phi_\sigma) \phi_H.
	\end{aligned}
$$
The parameters $\lambda_\alpha^\pro$ and $\lambda_\alpha^\apo$ are the proliferation and apoptosis rates corresponding to the $\alpha$-th species. Furthermore, $\lambda_{P\!H}$ denotes the transition rate from the proliferative to the hypoxic phase below the nutrient level $\sigma_{P\!H}$, $\lambda_{H\!P}$ the transition rate from the hypoxic to the proliferative phase above the nutrient level $\sigma_{H\!P}$, and $\lambda_{H\!N}$ the transition rate from the hypoxic to the necrotic phase below the nutrient level $\sigma_{H\!N}$. Lastly, $\calH$ represents the Heaviside step function that can be replaced by the Sigmoid function if a sufficiently smooth right-hand side is required.

In the instance of diffusion-type models, \cite{tao2011chemotaxis, tao2014energy}, \cite{engwer2017structured} and \cite{sfakianakis2020hybrid} discuss related theories of ECM degradation due to MDEs generated by hypoxic cell concentrations and subsequent tumor invasion and metastasis.  Following these references, we present an ECM evolution equation in the form of:
$$
	\boxed{\begin{gathered} 
			\textbf{Stratified tumor growth model with ECM: $\OD$} \\\begin{aligned}
	\pt \ecm ={}& S_{\ECM}( \phi_\A) \\
	={}&-\lambda_{\ECM}^\deg \ecm\mde + \lambda_{\ECM}^\pro \phi_\sigma (1- \ecm) \calH(\ecm - \phi_{\ECM}^\pro)
	\end{aligned}  \end{gathered} }
$$
Here, $\lambda_{\ECM}^\deg$ denotes the degradation rate of ECM fibers due to the matrix degrading enzymes, and $\lambda_{\ECM}^\pro$ is the production rate of ECM fibers above the threshold level $\phi_{\ECM}^\pro$. Further, for $(\phi_\beta)_{\beta \in \RD}$ we arrive at the following set of equations:
$$\boxed{\begin{gathered} 
	\textbf{Stratified tumor growth model with ECM: $\RD$} \\
	\begin{aligned}
	\pt \phi_\sigma + \div(\phi_\sigma v) ={}& \div \big(M_\sigma \nabla (D_\sigma \phi_\sigma - \chi_c  (\phi_P +\phi_H )\big)+S_\sigma( \phi_\A) \\
	\pt \mde                        ={}& M_{\MDE} D_{\MDE} \Delta \mde+S_{\MDE}( \phi_\A)
	\\
	\pt \taf                        ={}& M_{\TAF} D_{\TAF} \Delta \taf+S_{\TAF}( \phi_\A)                      
\end{aligned}\end{gathered}} $$
where the source functions are given by
$$
	\begin{aligned}
S_{\sigma}(\phi_\A) ={}&  \lambda_P^\apo \phi_P + \lambda_H^\apo \phi_H-\lambda_P^\pro \phi_\sigma\phi_P(1-\phi_T) -\lambda_{H}^\pro \phi_\sigma\phi_H(1-\phi_T)   \\
& + \lambda_{\ECM}^\deg \ecm\mde- \lambda_{\ECM}^\pro \phi_\sigma (1- \ecm) \calH(\ecm - \phi_{\ECM}^\pro) ,  \\
S_{\MDE}(\phi_\A) ={}&  \lambda_{\MDE}^\pro (\phi_P + \phi_H) \ecm\frac{\sigma_{H\!P}}{\sigma_{H\!P} + \phi_\sigma} (1-\mde)-\lambda_{\MDE}^\deg \mde  \\[-.25cm] &- \lambda_{\ECM}^\deg \ecm \mde,  \\
S_{\TAF}(\phi_\A) ={}& \lambda_{\TAF}^\pro (1- \taf) \phi_H \calH(\phi_H-\phi_{H}^\pro) - \lambda_{\TAF}^\deg \taf.  
\end{aligned}
$$
The parameters $\lambda_{\MDE}^\deg$ and $\lambda_{\TAF}^\deg$ denote the decay rates of the MDEs and TAFs, respectively. Moreover, $\lambda_{\MDE}^\pro$  represents the production rate of MDEs, and $\lambda_{\TAF}^\pro$ is the production rate of the $\taf$ due to the release by hypoxic cells above a threshold value of $\phi_{H}^\pro$. 

We notice that the cell species $\phi_\alpha$, $\alpha \in \{P,H,N,\sigma,\ECM\}$, form a mass conserving subsystem in the sense that their source terms add to zero. The constituents  $\phi_\MDE$ and $\phi_\TAF$ do not belong to a mass exchanging closed subsystem since they are signals and show natural degradation factors that are not absorbed by the other constituents.

\begin{figure}[htb!]
    \centering
    \includegraphics[width=.95\textwidth]{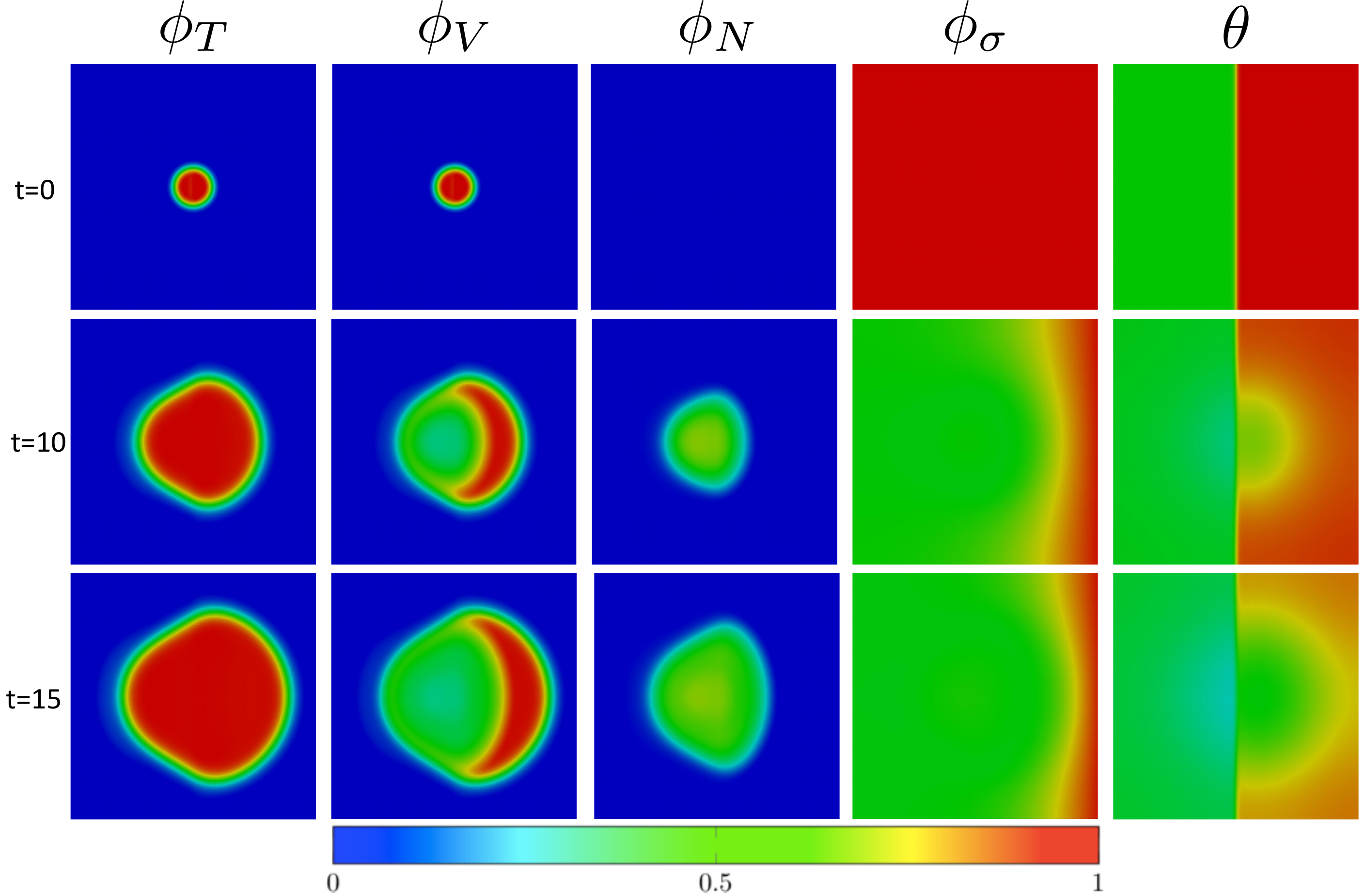}
    \caption{Evolution of the tumor mass under the influence of the extracellular matrix, taken with permission from Figure 1 and 2 in \cite{fritz2019unsteady}.}
    \label{Fig:Image_ECM}
\end{figure}

Numerically, we depict a simulation of a tumor with the degradation of the ECM in \cref{Fig:Image_ECM}. The viable part of the tumor consists of the proliferative and hypoxic phase. It absorbs the nutrients and starts to grow until $t=5$. Then the nutrients are sufficiently deprived in the sense that a necrotic core forms. The tumor moves towards the right and cell-to-cell and cell-to-matrix adhesion effects can be observed, i.e., tumor cells move towards nutrients due to chemotaxis and tumor cells move towards the ECM due to haptotaxis.

\subsection{Nonlocal phenomena} \label{Sec:Mod:Nonlocal}
In this section, the nonlocal impacts of tumor evolution models are discussed.
There are two types of nonlocality: spatial and temporal.
The first phenomenon relates to a time-fractional derivative in the PDE and is known as the memory effect.
In the second scenario, one must deal with a space integral, which reflects long-range interactions.

In addition, nonlocal events are incorporated into mathematical models of cancer cells.
These effects demonstrate long-distance interactions and may be geographical or temporal in character.
In the situation of spatial nonlocality, cell-to-matrix and cell-to-cell adhesion qualities are crucial to tumor growth modeling and encourage the proliferation of tumor cells.
Due to the structure of integro-differential systems, these events are nonlocal in space and require a special mathematical treatment.
\cite{fritz2019unsteady} explored cell-to-cell adhesion, whereas \cite{fritz2019local} investigated cell-to-matrix adhesion.
Further, we mention the articles by \cite{scarpa2021class} and \cite{frigeri2017diffuse} that studied nonlocal cell-to-cell adhesion properties  in phase-field models with applications to tumor growth.

In the case of temporal nonlocality, not only does the outcome of the previous step affect the current evolution, but it is also taken into account that cells have innate memories \citep{meir2020single}.
The past consequently effects the present.
In contrast to the normal Fickian diffusion process, memory effects are handled using a time-fractional derivative and fractional heat equations reflect the process of subdiffusivity.
As evidenced by  the in vitro and in vivo experimental findings of \cite{jiang2014anomalous}, tumors migrate via both traditional Fickian diffusion and subdiffusion.
\cite{fritz2021timefractional,fritz2021equivalence} investigated the memory effect in connection to the time-fractional Cahn--Hilliard equation with degenerating mobility.
Additionally, \cite{fritz2021subdiffusive} examined a fractional tumor model including subdiffusion, nutritional couplings, and mechanical deformation. 

\subsubsection{Nonlocal-in-space: cell-to-cell and cell-to-matrix adhesion}
If events or cell concentrations at one site in the tumor domain are dependent on events at other points within a defined neighborhood, it is said that the model is spatially nonlocal.
Long-distance interactions, such as cell-to-cell adhesion, are among the several processes that affect the mobility and migration of tumor cells.
cell-to-cell adhesion is a crucial aspect of tissue formation, stability, and degeneration, as well as a major contributor to cancer cell invasion and metastasis. 

Following \cite{chaplain2011mathematical} and \cite{frigeri2017diffuse}, we address cell-to-cell adhesion effects, which are responsible for the binding of two or more cells via protein processes on their respective cell surfaces.
The Ginzburg--Landau free energy functional generates separation and surface tension effects \citep{frigeri2017diffuse}, hence it is reasonable to incorporate cell-to-cell adhesion.
Therefore, tumor cells prefer to adhere to each other rather than healthy cells.
The physicists \cite{giacomin1996exact,giacomin1997phase} studied the problem of phase separation from a microscopic background using statistical mechanics and obtained the Helmholtz energy functional
$$\begin{aligned}
\mathcal{E}(\phi_T)={}&\int_\Omega \Psi(\phi_T)\, \text{d}x  + \frac14 \int_\Omega \int_\Omega J(x-y) \big(\phi_T(x)-\phi_T(y)\big)^2 \, \text{d}y \,\text{d}x.\end{aligned}$$
In this equation, we assume that $J:\mathbb{R}^d \to \mathbb{R}$ is a convolution kernel with the essential symmetry property $J(-x)=J(x)$. One obtains the Ginzburg--Landau energy by choosing a particular kernel sequence and passing the limit \citep{frigeri2015a}.
We modify the energy to account for chemotaxis and consider
\begin{equation*} \begin{aligned}
\mathcal{E}(\phi_T,\phi_\sigma) =&\int_\Omega \Psi(\phi_T) + \frac{D_\sigma}{2} \phi_\sigma^2 - \chi_c \phi_T \phi_\sigma \text{d}x \\ &+ \frac14 \int_\Omega \int_\Omega J(x-y) \big(\phi_T(x)-\phi_T(y)\big)^2 \, \text{d}y \, \text{d}x. \end{aligned}
\end{equation*}
Hence, we propose a class of long-range interactions, which are represented by chemical potentials of the form
\begin{equation*} \label{mu}
\mu_T = \frac{\delta \mathcal{E}}{\delta \phi_T} = \Psi'(\phi_T) - \chi_c \phi_\sigma + \int_\Omega J(x-y) \big(\phi_T(x)-\phi_T(y)\big) \, \textup{d} y.
\end{equation*}
This immediately results in the nonlocal system:
\begin{equation} \boxed{
	\begin{gathered} 
	\textbf{Four-species model with cell-to-cell adhesion} \\
	\begin{aligned}
	\partial_t \phi_T + \div(\phi_T v) ={}&  \div \big(M_T \phi_T^2(1-\phi_T)^2 \nabla \mu_T\big) + S_T(\phi_T,\phi_\sigma) \\
	\mu_T ={}& \Psi'(\phi_T)- \chi_c \phi_\sigma + \phi_T \cdot J*1 - J*\phi_T  \\
	\partial_t \phi_\sigma + \div(\phi_\sigma v) ={}&  \div\big(M_\sigma \nabla (D_\sigma  \phi_\sigma- \chi_c  \phi_T) \big) + S_\sigma(\phi_T,\phi_\sigma)
	\end{aligned}\end{gathered} }
\label{Eq:Nonlocal} \end{equation}

Included in models that account for cell-to-matrix adhesion effects are MDEs that erode the ECM; hence, this mechanism permits cell migration into tissue.
Such systems have been extensively investigated in \cite{engwer2017structured} and \cite{chaplain2011mathematical}.
In contrast to the fourth-order Cahn--Hilliard phase-field equation in our case, a reaction-diffusion equation is used to describe the tumor volume fraction in these publications.
The cell-to-matrix adhesion flux can be categorized as a local gradient-based haptotaxis effect \citep{stinner2014global,tao2011chemotaxis,walker2007global} or a nonlocal adhesion-based haptotaxis effect \citep{armstrong2006continuum,chaplain2011mathematical,gerisch2008mathematical}.
We consider the respective fluxes of the form
$$
	J_\alpha(\phi_\A) = \chi_h \phi_V \cdot \begin{cases} \nabla \ecm, &\alpha = \text{local}, \\ k*\ecm, &\alpha=\text{nonlocal},\end{cases}
$$
where $k$ is a vector-valued kernel function. This adhesion flux is factored into the equation of the extended mass balance law for the volume fraction of viable cells as follows: 
$$\boxed{\pt \phi_V + \div(\phi_V v) = \div(m_V(\phi_\A) \nabla \mu_V) + \div J_\alpha + S_V(\phi_\A)}$$

\begin{figure}[htb!]
    \centering
    \includegraphics[width=.9\textwidth,page=4]{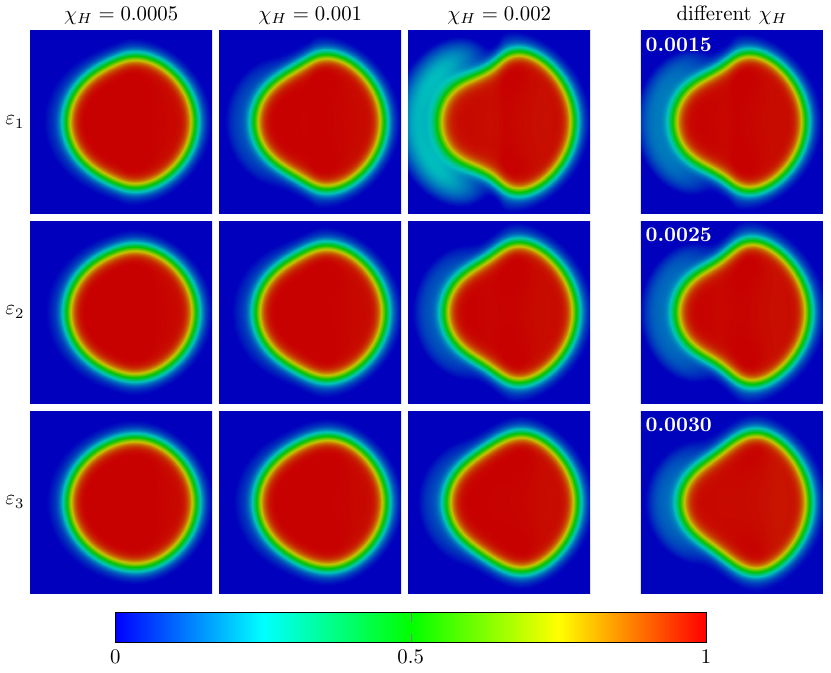}
    \caption{Simulation of the tumor volume fraction $\phi_T$ for different haptotaxis parameters $\chi_h \in \{0.0005, 0.001, 0.002\}$ and different values of $\eps \in \{\eps_1,\eps_2,\eps_3\}=\{0,0.00275,0.00525\}$ for a fixed time; taken with permission from Figure 5 in \cite{fritz2019local}.}
    \label{Fig:Image_Nonlocal}
\end{figure}

In the following, we numerically investigate the effects of the different haptotaxis parameters on the growth of the tumor mass in \cref{Fig:Image_Nonlocal}. We distinguish between three different values, and we compare the local gradient-based adhesion flux to the nonlocal one with two different values of $\eps$. We notice that different pairing of $\eps$ and $\chi_h$ result in similar simulations. A larger value of $\chi_h$ results in a larger difference between the local and nonlocal model. 
Comparing the local and nonlocal models, we find a greater adhesion impact in the local model, as the tumor mass moves further to the right of the boundary where nourishment is put in the local model. By selecting a smaller haptotaxis parameter, the local model can resemble the nonlocal model for a fixed value of $\eps$. The subtleties of the discretization of the nonlocal model is further discussed in \cref{Sec:Numerics}. 

\subsubsection{Nonlocal-in-time: The memory effect}
According to \cite{Balkwill12}, \cite{wang17} and \cite{yuan16}, the tumor microenvironment significantly influences the proliferation and migration of tumor cells.
In addition to Fickian diffusion and subdiffusion, tumor cells travel through a variety of methods.
The results of the experiments of \cite{jiang2014anomalous} show anomalous diffusion in the progression of cancer.
In addition to clinical data from patients with adrenal and liver tumors, they discovered subdiffusion during in vitro tests of generating cultured cells from the breast line and during in vitro trials of developing cultured cells from the liver line.

In earlier sections, the phenomenological law $J_T=-m_T(\phi_\A)\nabla\mu_T$ was used to depict the typical relationship between flow and the gradient of the chemical potential.
A more complicated phenomenological link that could account for hypothesized nonlocal, nonlinear, and memory effects \citep{gorenflo2002time, povstenko2017two}, can be substituted for this law without contradicting the conservation law suggested by the continuity equation.
\cite{seki2003recombination} and \cite{yuste2004reaction} simulate subdiffusion-limited reactions on a tiny scale by employing fractional derivatives in flux and reaction terms.
Consequently, we propose
$$
\begin{aligned}
J^\rel_T(\phi_\A)&=-\p_t \big (g_\alpha \circledast (m_T(\phi_\A)\nabla\mu_T)\big),\\
S^\rel_T(\phi_\A)&=\p_t \big(g_\alpha \circledast S_T(\phi_\A) \big),
\end{aligned}
$$
for $\alpha\in(0,1)$. Inserting the relaxed flux and source  into the law of conservation of mass \cref{Eq:MassBalance}, it yields for the tumor species $\phi_T$ 
\begin{equation*}
\begin{aligned}
\p_t \phi_T ={}& -\div J_T^\rel(\phi_\A) + S_T^\rel(\phi_\A) = \p_t \big( g_\alpha \circledast \big( \div(m_T(\phi_\A)\nabla\mu_T) +  S_T(\phi_\A) \big) \big).
\end{aligned}
\end{equation*}
We can equivalently rewrite this equation by taking the convolution with $g_{1-\alpha}$ on both sides of the equation and using the inverse convolution property \cref{Eq:InverseConvolution}. Therefore, we obtain
\begin{equation*}  \label{Eq:CHCaputo}
\boxed{\begin{aligned}
\p_t^\alpha \phi_T ={}& \div(m_T(\phi_\A)\nabla\mu_T) + S_T(\phi_\A).
\end{aligned}}
\end{equation*}
The chemical potential reads $\mu_T=\Psi'(\phi_T)-\eps_T^2\Delta \phi_T$ in the case of the Ginzburg--Landau energy \cref{Eq:Ginzburg}. This model is called the time-fractional Cahn--Hilliard equation  \citep{fritz2021timefractional,fritz2021equivalence}. If one selects the Dirichlet energy $\calE(\phi_T)=\int_\Omega \phi_T^2 \dx$, then one obtains the time-fractional reaction-diffusion equation as studied in a tumor growth setting in \cite{fritz2021subdiffusive}.

Typically, subdiffusive models do not exhibit linear growth if enough nutrients are in the tumor environment, as it can be typically observed for integer-order tumor growth models. Indeed, subdiffusive models have a larger growth in the beginning and the growth becomes damped afterwards. This can be explained by the memory effect of cells that first try to absorb everything that they can get and, afterwards, become more lenient with the available nutrients. The details for the implementation of the time-fractional derivative are considered in \cref{Sec:Numerics}.
\begin{figure}[htb!]
    \centering
    \includegraphics[width=.6\textwidth]{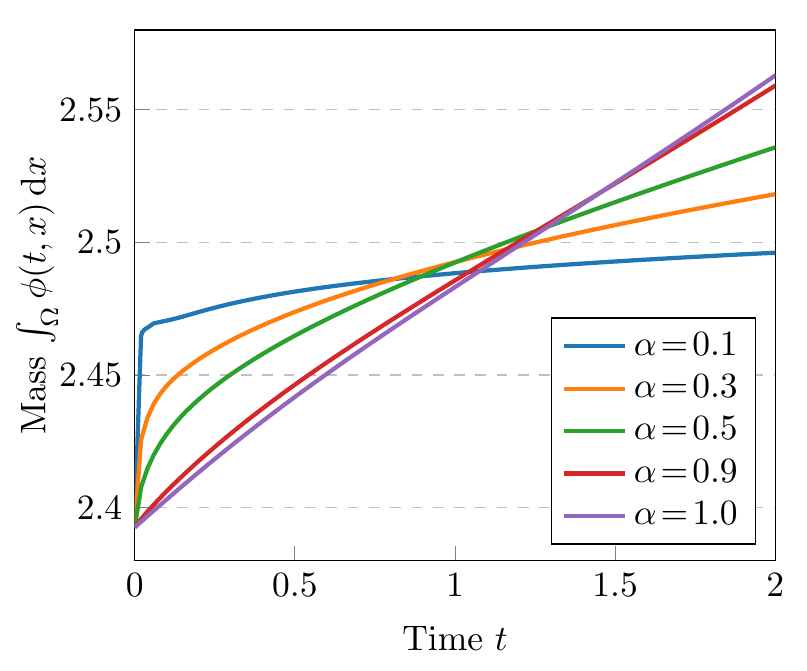}
    \caption{Evolution of tumor mass $\phi_T$ with different parameters of the fractional order $\alpha$; taken with permission from Figure 5 in \cite{fritz2021timefractional}.}
    \label{Fig:Image_Frac}
\end{figure}

\subsection{Uncertainty in tumor modeling} \label{Sec:Mod:Stochastic}
In thin subdomains at the interfaces of the phase fields, stochastic variations of the phase concentrations are possible, see the works by \cite{orrieri2020optimal} and \cite{fritz2021modeling,fritz2023stochastic}. The variances of these regions of random behavior are constrained by noise parameters $\phi_\alpha^\omega$ with noise intensity $\omega_\alpha$. In fact, the variations in $\phi_\alpha$ are restricted to interface regions by making use of the operator
$$G_\alpha(\phi_\alpha)=\omega_\alpha \mathcal{H}((\phi_\alpha-\phi_\alpha^\omega)(1-\phi_\alpha-\phi_\alpha^\omega)).$$ %
Generally, the model incorporates the randomness in the evolution of species along the interface as a cylindrical Wiener process on $L^2(\Omega)$. We refer to \cite{da1996stochastic}, \cite{cardon2001cahn}, and \cite{elezovic1991stochastic} on the stochastic Cahn--Hilliard equation. Further details on stochastic PDEs can be found in the textbooks \cite{lord2014introduction}, \cite{prevot2007concise}, and \cite{liu2015stochastic}. 

Modeling-wise, we add $G_\alpha \dot W_\alpha$  to the mass balance equation for $\phi_\alpha$ and to keep the mass balance equations in standard form, we slightly abuse the standard notation by writing $\dot W_\alpha$ in the sense $\dot W_\alpha \dt = \dd W_\alpha$. In the case of the simplified four-species model, we obtain the following stochastic version of it.

$$\label{Eq:Wiener} \boxed{
		\begin{gathered}
	\textbf{Stochastic four-species model} \\
	\begin{aligned} 
	\pt \phi_T \!+\! \div(\phi_T v) ={}& \div (M_T \phi_T^2(1-\phi_T)^2 \nabla \mu_T) + S_T( \phi_T,\phi_\sigma) + G_\alpha(\phi_T) \omega_T \text{d}W_T    \\
	\mu_T ={}& \Psi'(\phi_T) - \chi_c \phi_\sigma- \eps_T^2 \Delta \phi_T \\
	\pt \phi_\sigma \!+\!  \div(\phi_\sigma v)                     ={}& \div\Big(M_\sigma \nabla\big( D_{\sigma} \phi_\sigma -\chi_c\phi_T  \big) \Big) +S_{\sigma}( \phi_T,\phi_\sigma) \\
	\end{aligned}\end{gathered} } $$

\begin{figure}[htb!]
    \centering
    \includegraphics[width=.9\textwidth]{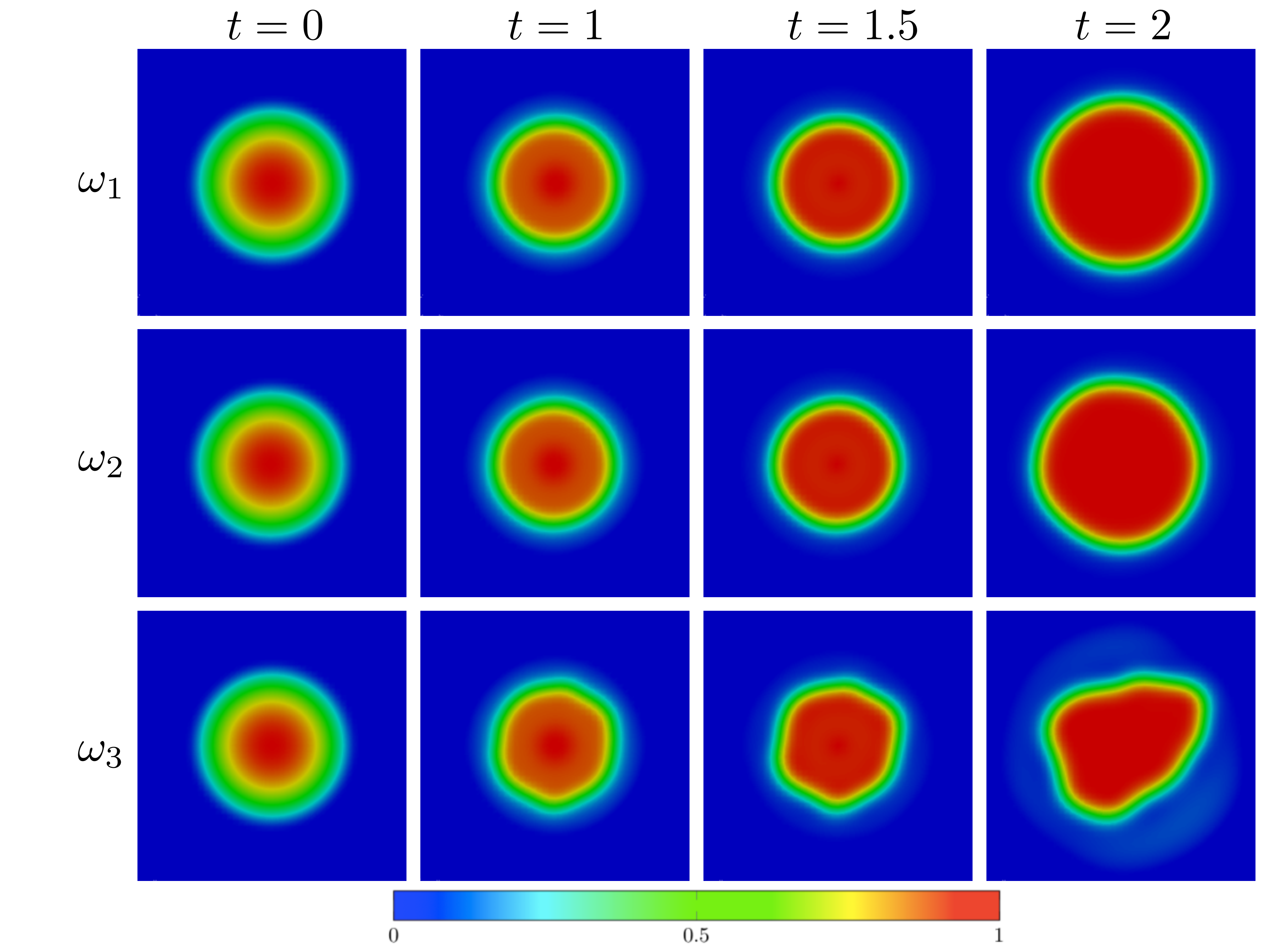}
    \caption{Evolution of the tumor volume fraction for different values of the noise intensity $\omega_T$; we consider the cases of $\omega_T\in\{\omega_1,\omega_2,\omega_3\}=\{0,0.001,0.1\}$.}
    \label{Fig:stochastic}
\end{figure}

 Numerically, we investigate the influence of the stochasticity in the tumor growth model in \cref{Fig:stochastic}. We first consider the deterministic model that corresponds to $\omega_T=0$ and, afterwards, compare it to the cases with two different values for $\omega_T$. We notice that a larger value of $\omega_T$ results in a non-regular shape of the tumor's interface. The implementation of the Wiener process is discussed in \cref{Sec:Numerics}.

\subsection{Mechanical deformation} \label{Sec:Mod:Mechanical}
As the tumor grows, the surrounding host tissues generate mechanical stress, restricting the tumor's growth.
In the papers \citep{faghihi2020coupled,lima2016selection,lima2017selection}, mechanical deformation in a tumor development model was first mentioned, and in terms of analysis, it was first examined in \cite{fritz2021subdiffusive} in a diffusion-type tumor model and subsequently \cite{garcke2021phase} in a Cahn--Hilliard type system.
Such models with elasticity are referred to as Cahn--Larch\'e equations. It had previously been incorporated into the Cahn--Hilliard equation without being applied to tumor growth or traditional source variables in \cite{garcke2003cahn,garcke2005cahn}.

As the tumor grows, the surrounding host tissues generate mechanical stress, restricting the tumor's growth.
Regarding mathematical modeling and sensitivity studies, several works \citep{lima2016selection,lima2017selection,hormuth2018mechanically,faghihi2020coupled} have employed reaction-diffusion equations with mechanical coupling to predict tumor progression.
\cite{fritz2021subdiffusive} examined the well-posedness of a model in which similar mechanical factors were incorporated.
The underlying energy functional now contains the stored energy potential $W(\phi_T,\eps(u))$, which is dependent on the tumor volume fraction $\phi_T$ and the symmetric strain measure $\eps(u)=\frac12(\nabla u+\nabla u^\top)$ of the displacement field $u$.
Assuming minor deformations, we consider the specific stored energy potential 
\begin{equation} \label{Eq:Stored}
W(\phi_T,\eps(u)) = \frac{1}{2} \eps(u):T_M(\phi_T)\eps(u) + \eps(u):T_S(\phi_T),
\end{equation}
where $T_S(\phi_T)=\lambda \phi_T \bbI$ is the symmetric compositional stress tensor with $\lambda>0$ and $T_M$ is the linear elastic inhomogeneous material tensor.
The symbol $\bbI$ signifies the $(d\times d)$-dimensional identity matrix in this instance.
The displacement field $u$ is governed by the conservation equation of linear and angular momentum
$$
\begin{aligned}
\p_t (\phi_T v) + \div(\phi_T v\otimes v) &= \div\, T_C + \phi_T b + p, \\
T_C- T_C^\top &= \text{m},
\end{aligned}
$$
where $v$ is the volume-averaged velocity, $b$ is the body force, $p$ is the momentum contributed by other components, and $\text{m}$ is the intrinsic moment of momentum.
First variations of the energy functional $\mathcal{E}$ with respect to $\phi_T$ and $\eps(u)$, respectively, determine the chemical potential $\mu_T$ and Cauchy stress tensor $T_C$.
We minimize the system's complexity by using the typical simplification assumptions of \cite{lima2016selection}.
Particularly, we assume constant mass density $\text{m}=0$ and a monopolar material $b=0$.
In addition, we disregard inertial forces and set $\div(\phi_T v\otimes v)=p=0$.
We assume that the mechanical equilibrium is reached quicker than diffusion, i.e., that the time derivative on the left-hand side disappears.
After the simplifications, the mechanical deformation equation \cref{Eq:Stored} becomes 
$$0 = \div\, T_C = \div\, \frac{\delta \calE(\phi_T,\phi_\sigma,\eps(u))}{\delta \eps(u)} = \div\, \frac{\partial W(\phi_T,\eps(u))}{\partial \eps(u)}.$$
We assume that the tumor is an isotropic and homogeneous material, i.e., that its material tensor $C_M(\phi)=C_M$ has the form $$C_M\eps(u)=2G\eps(u)+\frac{2G\nu}{1-2\nu} \tr\,\eps(u)\bbI,$$ where $G>0$ and $\nu<\frac{1}{2}$ represent the shear modulus and Poisson ratio, respectively.
For the stored energy potential $$
W(\phi_T,\eps(u)) = \frac{1}{2} \eps(u):(2G\eps(u)+\frac{2G\nu}{1-2\nu} \tr\,\eps(u)\bbI)\eps(u) + \eps(u):(\lambda \phi_T \bbI).
$$
and its partial derivatives with respect to $\phi_T$ and $\eps(u)$, we may therefore write 
$$\begin{aligned}\frac{\partial W(\phi_T,\eps(u))}{\partial \phi_T} &= \lambda \div u, \\
\frac{\p W(\phi_T,\eps(u))}{\p \eps(u)} &= 2 G \eps(u) + \frac{2 G \nu}{1-2\nu} \tr(\eps(u)) \bbI + \lambda \phi_T \bbI.
\end{aligned}$$

With mechanical deformation, it provides the model
$$  \boxed{
	\begin{gathered}
	\textbf{Four-species model with mechanical deformation} \\
	 \begin{aligned}
	\pt \phi_T+ \div(\phi_T v) ={}& \div \big(M_T \phi_T^2(1-\phi_T)^2  \nabla \mu_T\big) + S_T( \phi_T,\phi_\sigma)    \\
	\pt \phi_\sigma +  \div(\phi_\sigma v)                     ={}& \div\big(M_\sigma \nabla( D_{\sigma} \phi_\sigma -\chi_c\phi_T  ) \big) +S_{\sigma}( \phi_T,\phi_\sigma) \\
	0={}& \div \Big(2 G \eps(u) + \frac{2 G \nu}{1-2\nu} \tr(\eps(u)) \bbI + \lambda \phi_T \bbI\Big)  
	\end{aligned} \end{gathered}} $$
As \cite{fritz2021subdiffusive} demonstrated in their study, the Ginzburg--Landau energy yields $$\mu_T = \Psi'(\phi_T) - \chi_c \phi_\sigma- \eps_T^2 \Delta \phi_T + \lambda \div u,$$ whereas the Dirichlet energy yields $\mu_T =D_T \phi_T - \chi_c \phi_\sigma + \lambda \div u$, i.e.,
$$\calE(\phi_T,\phi_\sigma)= \int_\Omega \Big\{ \frac{D_T}{2} \phi_T^2 + \frac{D_\sigma}{2} \phi_\sigma^2 - \chi_c \phi_T \phi_\sigma + W(\phi_T,\eps(u))\Big\} \dx.$$

\subsection{Chemotherapeutic influence} \label{Sec:Mod:Chemo}
In addition to precisely simulating the tumor's growth, mathematicians are interested in treating the tumor and stopping its growth.
Currently, chemotherapy, surgery, immunotherapy, and radiotherapy are used to treat malignancies.
Angiogenesis is one of the primary mechanisms by which tumors grow, hence anti-angiogenic drugs that inhibit the production of new vascular structures are commonly identified as one of the methods to delay or stop cancer growth.
Consequently, a realistic model of angiogenesis is essential for evaluating the efficiency of anti-angiogenic drugs; for the ideal dosage of medication, see the optimal control problems discussed in \cite{colli2020mathematical,colli2021optimal}.
Chemotherapy was incorporated into our research \citep{fritz2021subdiffusive} with a reaction-diffusion equation and subdiffusive tumor growth, as well as in the articles \citep{ebenbeck2019optimal,signori2019penalisation,garcke2018optimal,colli2020mathematical,colli2021optimal} on optimum control issues for the optimal drug dosage. Moreover, in the work \cite{wagner2023phase} it was assumed that the immunotherapeutic concentration follows the Hill–Langmuir equation. 

In addition to studying the growth of tumors, we also incorporate a substance that inhibits their spread.
Current cancer treatments include: \medskip 
\begin{itemize}  \itemsep0.3em
	\item Surgery: Removing the tumor by an operation. 
	\item Immunotherapy: Strengthening the immune system.
	\item Radiotherapy: Employing radiation to eradicate cancerous cells. 
	\item Chemotherapy: Utilizing medications to destroy the tumor. \medskip
\end{itemize} 
These therapies, excepting surgery, are administered in cycles, with each cycle consisting of a period of therapy followed by a period of rest to allow the patient's body to repair and regenerate new, healthy cells.
These therapeutic procedures should diminish the tumor to a degree where surgical removal is feasible. 

The mass density of chemotherapy $\cmt$ is considered to be driven by a reaction-diffusion equation that links to the tumor equation and, if chemotherapy is present, degrades the tumor.
Therefore, we recommend adding the index $\CMT$ to the index set $\RD$ and propose the model:
$$\label{Eq:Chemo} \boxed{
		\begin{gathered}
	\textbf{Four-species model with chemotherapy} \\
	\begin{aligned} 
	\pt \phi_T+ \div(\phi_T v) ={}& \div (M_T \phi_T^2(1-\phi_T)^2 \nabla \mu_T) + S_T( \phi_T,\phi_\sigma,\cmt)    \\
	\mu_T ={}& \Psi'(\phi_T) - \chi_c \phi_\sigma- \eps_T^2 \Delta \phi_T \\
	\pt \phi_\sigma +  \div(\phi_\sigma v)                     ={}& \div\Big(M_\sigma \nabla\big( D_{\sigma} \phi_\sigma -\chi_c\phi_T  \big) \Big) +S_{\sigma}( \phi_T,\phi_\sigma,\cmt) \\
	\pt \cmt ={}& M_\CMT D_\CMT  \Delta \cmt + S_\CMT(\phi_T,\phi_\sigma,\cmt).
	\end{aligned}\end{gathered} } $$
The mobility of chemotherapeutic agents  is given by $M_\CMT$ and the source $S_\CMT$ reads
$$S_\CMT(\phi_T,\phi_\sigma,\cmt)= -\lambda_\CMT^\deg \cmt - \lambda_\CMT^\kill \frac{\phi_T (1-\phi_T) \cmt}{K_\CMT+\cmt},$$
where $\lambda_\CMT^\deg$ is the degradation factor of chemotherapeutic agents and $\lambda_\CMT^\kill$ represents the rate at which chemotherapeutic agents act and are subsequently blocked by the death of tumor cells. 
The killing term includes a saturation effect, so that chemotherapy is most effective against cells in a certain growth phase.
The parameter $K_\CMT>0$ is the density of chemotherapeutic drugs at half-maximum concentration.
Similarly, the source term of the tumor volume fraction will contain a term of the kind $$-\lambda_T^\kill \frac{\phi_T (1-\phi_T) \cmt}{K_\CMT+\cmt},$$ that represents the chemotherapy's killing impact at some rate $\lambda_T^\kill$.
In our approach in \cite{fritz2021subdiffusive}, chemotherapeutic agents in cycles are represented by a Dirichlet border of the type X. \\[-1cm]

$$\cmt(t,x)\vert_{x \in \p\Omega}=\begin{cases} 1, &\text{for $t\leq 2$ or $6<t \leq 8$ or $12<t\leq 14$,} \\ 0, &\text{else}.\end{cases}$$
That is, during the time $t \in [0,2] \cup (6,8] \cup (12,14]$ chemotherapy treatment is provided and in between, the body is permitted to rest.

\subsection{Angiogenesis and mixed-dimensional coupling} \label{Mod:Mixed}
Hypoxic tumor cells not only release MDEs to degrade the ECM, but also TAFs, which stimulate endothelial cell proliferation and new vessel formation.
Angiogenesis is the process of blood vessels sprouting and elongating in order to supply the tumor with nutrients.
Unless sufficient nutrients and oxygen are available for proliferation, the volume of an isolated colony of tumor cells is often limited to 1\si{mm}$^3$, as shown in \cite{nishida2006angiogenesis}, unless adequate nutrients and oxygen are supplied.
In order to access these nutrients, cancerous cells drive angiogenesis \citep{carmeliet2011molecular, patsch2015generation}.
Regarding angiogenesis modeling and numerical simulations, we refer to \cite{cristini2009nonlinear,cristini2010multiscale,xu2016mathematical}.
In \cite{fritz2021analysis}, we studied angiogenesis in terms of the mathematical analysis of weak solutions in a Cahn--Hilliard-type model.
We are unaware of any subsequent works.
Due to mixed-dimensional couplings and the presence of hypoxic tumor cells that generate TAFs, these models are extremely complex. 

\cite{lima2014hybrid, xu2016mathematical, xu2017full, xu2020phase, wise2008three, santagiuliana2016simulation, santagiuliana2019coupling} presents the effect of angiogenesis on models of stratified tumor development.
In contrast to their prior techniques employing, for example, agent-based systems, we represent the network of blood arteries feeding a solid tumor mass as a network of 1D capillaries within a 3D tissue domain in our studies in \cite{fritz2021analysis,fritz2021modeling}.
In this perspective, tumor growth is viewed as a phase-field system with multiple cell species and other components.
The microvascular network in tumor-bearing tissue is modeled as a graph with 1D filaments through which nutrient-rich blood can flow.
This microvascular network is represented by $\Lambda$ and the individual edges by $\Lambda_i$, such that $\Lambda$ is given by the union $
\Lambda = \bigcup_{i=1}^N \Lambda_i.
$
An edge $\Lambda_i$ is parameterized with a curve parameter $s_i$ as follows:
$$
\Lambda_i = \left\{ x \in \Omega : x = \Lambda_i ( s_i ) = x_{i,1} + s_i \cdot ( x_{i,2}-x_{i,1} ),\;s_i \in (0,1 ) \right\}.
$$

\begin{figure}[htb!]
    \centering
    \includegraphics[width=.24\textwidth]{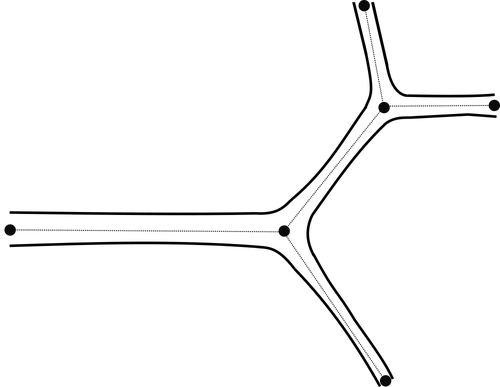}
    \includegraphics[width=.24\textwidth]{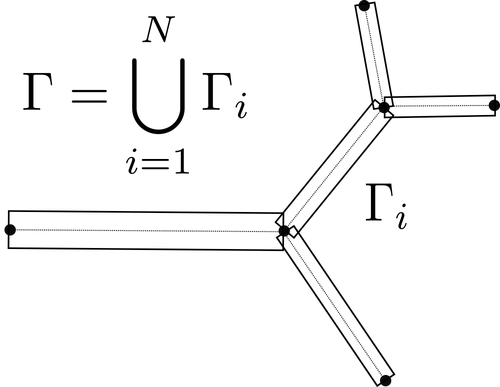}
    \includegraphics[width=.24\textwidth]{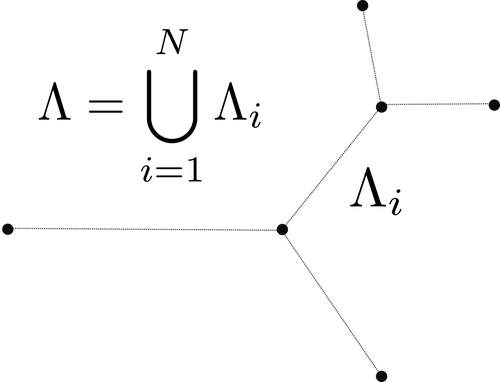}
    \includegraphics[width=.24\textwidth]{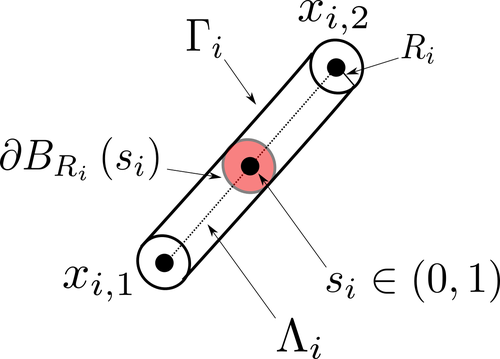}
    \caption{Discretization of the vessels into $N$-many vessels and introduction of 1D lines $\Lambda_i$; further, the vessel surfaces $\Gamma_i$ are depicted; taken with permission from Figure 2 in \cite{fritz2021analysis}.}
    \label{Fig:3d1d}
\end{figure}

We suggest $s$ as the global curve parameter for the entire 1D network $\Lambda$ by setting $s=s_i$ if $x = \Lambda(s) = \Lambda_i(s_i)$.
We search the domain $\Omega$ for 1D items that couple to their 3D counterparts for each value of the curve parameter $s$.
We suppose that the surface of a single vessel is a cylinder with a constant radius, and that the radius of a vessel attached to the edge $\Lambda_i$ is $R_i$.
We describe $\Gamma_i$ as the surface of the cylinder, with the edge $\Lambda_i$ as its center line, and the total surface $\Gamma$ is the union of the surfaces of the individual vessels $\Gamma_i$, see as well \cref{Fig:3d1d} for a depiction of the individual fields.

On the 1D network $\Lambda$, we take the constituents $\phi_v$, $v_v$ and $p_v$ into account, which reflect the 1D equivalents of the local nutrient concentration $\phi_\sigma$, the volume-averaged velocity $v$, and the pressure $p$.
We incorporate a new source term $S_{\sigma v}$ for coupling the 1D constituents $\phi_v$ and $p_v$ into the $\phi_\sigma$-equation.
Consequently, this source word is accountable for the relationship between the elements in $\Omega$ and $\Lambda$. 

To quantify the flux of nutrients across the vessel surface, we employ the Kedem--Katchalsky law \citep{ginzburg1963frictional} and write the flux $J_{\sigma v}$ between the nutrients on the network and tissue as \begin{equation} \label{Eq:Source3d1d}
J_{\sigma v} (\ov \phi_\sigma, \ov p, \phi_v, p_v ) = ( 1-r_{\sigma} ) f(\phi_{\sigma},\phi_v) L_p ( p_v - \overline{p} )  + L_\sigma ( \phi_v - \ov \phi_\sigma ),
\end{equation} 
where $r_\sigma>0$ is the reflection parameter, $L_\sigma,L_p>0$ represent the permeabilities of the vessel wall, and the function $f$ is either $\phi_\sigma$ or $\phi_v$ depending on the values of $p$ and $p_v$.
In addition, $\overline{p}$ represents the average circumferential pressure of cylinder cross-sections.
The averaging reflects the fact that the 3D-1D coupling is a reduced model from a physical standpoint, whereas the exchange occurs through the surface in a fully linked 3D-3D model.
The first portion of the Kedem--Katchalsky law measures the nutritional flux caused by the passage of blood plasma from arteries to tissues or vice versa.
It is defined by Starling's law, which is given by the pressure difference between $p_v$ and $p$ multiplied by a parameter $L_p$ representing the permeability of the vessel wall.
The second component of the law is a Fickian-type law that accounts for the tendency of nutrient concentrations to equalize. 

As the exchange activities between the vascular network and the tissue occur at the vessel surface $\Gamma$, we concentrate the flux $J_{\sigma v}$ using the Dirac measure $\delta_\Gamma$, i.e., by defining 
$$\langle \delta_\Gamma, \varphi \rangle_{C_c^\infty(\Omega)} = \int_\Gamma \varphi\vert_{\Gamma}(x) \, \dd S \quad \forall \varphi \in C_c^\infty(\Omega),$$
where $(C_c^\infty(\Omega))'$ is the space of distributions.
The resulting new source term in the nutrient equation is as follows: 
$$
S_{\sigma v}(\phi_\sigma,p,\phi_v,p_v) = J_{\sigma v} ( \phi_\sigma, p, \Pi_\Gamma \phi_v, \Pi_\Gamma p_v) \delta_\Gamma,
$$
where $\Pi_\Gamma \in \scrL(L^2(\Lambda);L^2(\Gamma))$ is the projection of the 1D quantities onto the cylindrical surface $\Gamma$ by extending the function value $\Pi_\Gamma \phi_v(s) = \phi_v(s_i)$ for all $s \in \p B_{R_i}(s_i)$. 

The 3D model reads:
$$  \boxed{\begin{gathered} 
\textbf{Angiogenesis model: 3D} \\
\begin{aligned}
\pt \phi_\alpha+ \div(\phi_\alpha v)            ={}& \div \big(m_\alpha( \phi_\A) \nabla \mu_\alpha\big)+ S_\alpha( \phi_\A)                                                      \\
\mu_\alpha                                 ={}&   \p_{\phi_\alpha} \Psi(\phi_\CH) - \eps^2_\alpha \Delta \phi_\alpha - \chi_c \phi_\sigma-\chi_h \ecm                                                      \\
\pt \phi_\beta                            ={}& S_\beta( \phi_\A)                                                                                                                                   \\
\pt \phi_\gamma    ={}& \div \big(m_{\gamma}( \phi_\A) D_{\gamma} \nabla \phi_\gamma\big)+S_{\gamma}( \phi_\A) \\
\pt \phi_\sigma + \div(\phi_\sigma v) ={}& \div \big(m_\sigma( \phi_\A) \nabla (D_\sigma  \phi_\sigma - \chi_c (\phi_P +\phi_H  )\big)+S_\sigma( \phi_\A)  \\ &+ J_{\sigma v} (\tr_\Gamma \phi_\sigma, \tr_\Gamma p, \Pi_\Gamma \phi_v,  \Pi_\Gamma p_v ) \delta_\Gamma \\
v  ={}& - K\big(\nabla p -S_p(\phi_\A,\mu_P,\mu_H)\big)  \\
\div\, v ={}&    L_p (\Pi_\Gamma p_v - p) \delta_\Gamma%
\end{aligned}\end{gathered}} $$
for $\alpha \in \{P,H\}$, $\beta \in \{N,\ECM\}$, $\gamma \in \{\MDE,\TAF\}$.

Since the vascular network is often composed of small inclusions, we averaged all the physical units across the cross-sections of the individual blood vessels and held them constant with regard to the angular and radial components.
In other words, the 1D variables $\phi_v$ and $p_v$ of a 1D vessel $\Lambda_i$ are entirely dependent on $s_i$.
\cite{koppl20203d} contains further information regarding the derivation of 1D pipe flow and transport models.
Consequently, the 1D model equations for vessel flow and transport are as follows: 
$$  \label{Eq:Model1D}
\boxed{\begin{gathered} 
	\textbf{Angiogenesis model: 1D} \\\begin{aligned}
\pt \phi_v + \p_{s_i} (v_v \phi_v) ={}&  \p_{s_i} (m_v(\phi_v)D_v \p_{s_i} \phi_v) -2\pi R_i J_{\sigma v} (\ov \phi_\sigma, \ov p, \phi_v,  p_v)\\
-  \;\p_{s_i} ( R_i^2 \pi K_{v,i} \; \p_{s_i} p_v )   ={}& -2 \pi R_i J_{pv}( \overline{p}, p_v ) \\
v_v ={}& - R_i^2 \pi K_{v,i} \p_{s_i} p_v
\end{aligned}\end{gathered}}
$$
In order to interconnect the multiple solutions on the vessel $\Lambda_i$ at inner network nodes at junctions $x \in \p \Lambda_i \setminus \p \Lambda$, we require the continuity of pressure and concentration as well as the conservation of mass, as shown in \cite{fritz2021modeling}.

\begin{figure}[htb!]
    \centering
    \includegraphics[width=.27\textwidth]{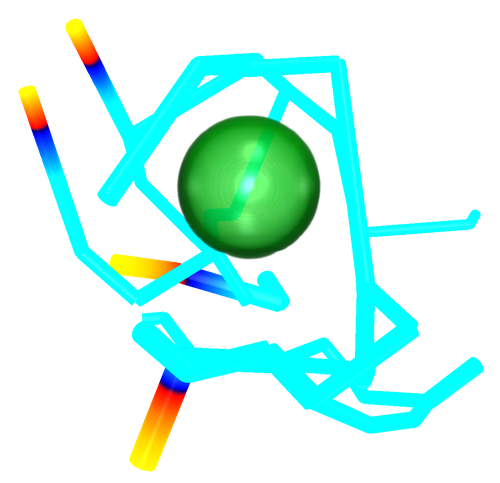} \!\!\!\!\!\!\!\!
    \includegraphics[width=.27\textwidth]{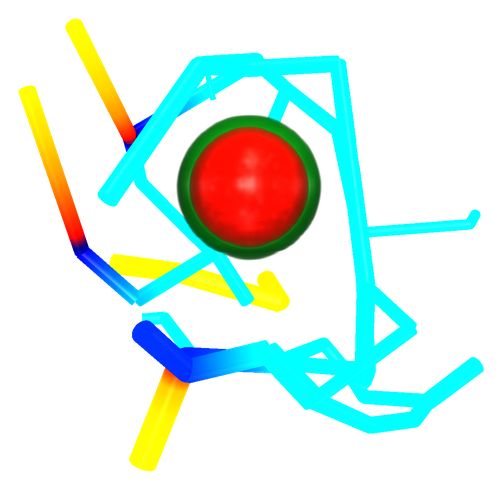}\!\!\!\!\!\!\!\!
    \includegraphics[width=.27\textwidth]{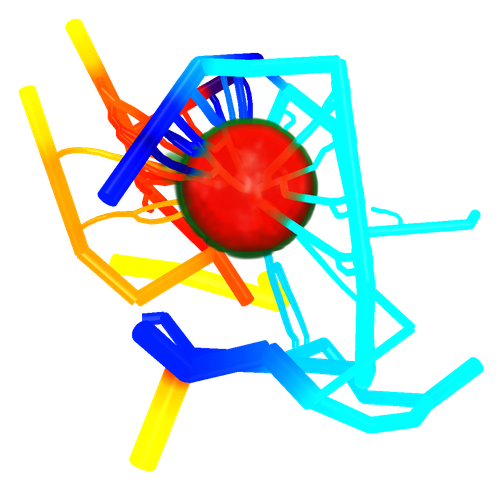}\!\!\!\!\!\!\!\!
    \includegraphics[width=.27\textwidth]{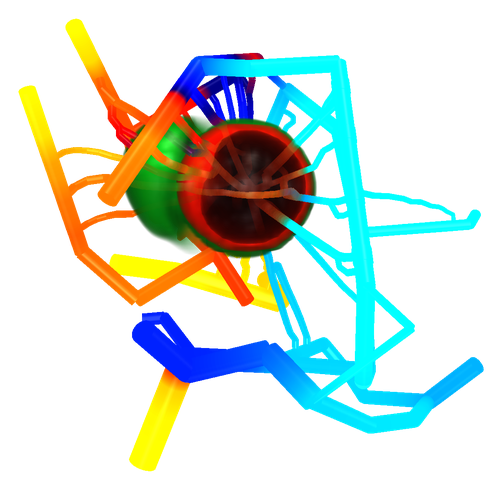} \\
    \includegraphics[width=.4\textwidth]{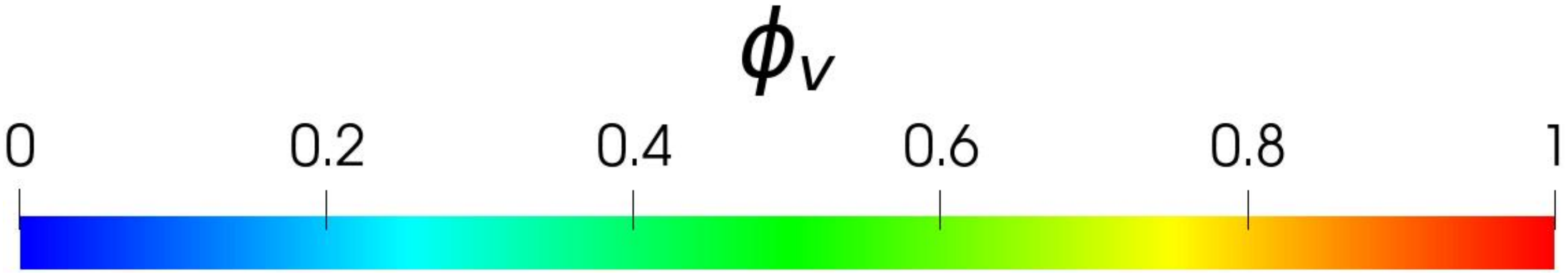}
    \caption{Stratification of tumor cells in its proliferative (green), hypoxic (red) and necrotic (black) phases over the time steps $t \in \{0,5,10,15\}$; growth of capillaries and movement of tumor cells to high-nutrient regions on the 1D lines that is expressed by the nutrients $\phi_v$.}
    \label{Fig:Angiogenesis}
\end{figure}

We present a numerical simulation of the tumor evolution in the setting of a capillary network in \cref{Fig:Angiogenesis}. We can observe that the tumor is deprived of nutrients and therefore, it stratifies, and it becomes hypoxic. The hypoxic tumor phase releases TAFs, and it is clear to see that angiogenesis happens -- the capillary begins to grow towards the tumor and provides it with new nutrients living in the 1D vessel network.

    	\section{Numerical Implementation} \label{Sec:Numerics}

\noindent Besides the analytical methods in the next section, we are interested in showing numerical simulations and studying the influence of the new features and effects of the models as we have seen in the previous section. How useful is a well-posed model that does not reflect real biological processes? In this section, we briefly describe the techniques that we previously used for the implementation of the PDEs in the last sections. Our code is based on the finite element libraries \texttt{libMesh} \citep{libMeshPaper} and \texttt{FEniCS} \citep{fenics}.  \texttt{FEniCS}  is written in the accessible \texttt{Python} language and variational forms are straightforward to implement. However, \texttt{libMesh} is a high
performance computing (HPC) library written in \texttt{C++} and therefore, yields higher potential for code optimization and saving run times than in \texttt{FEniCS}. We refer to our \texttt{GitHub} \medskip

\begin{center}  \url{https://github.com/CancerModeling/Angiogenesis3D1D} \medskip \end{center} 

\noindent where the code is freely accessible. In particular, the settings for the simulations in \cite{fritz2021analysis,fritz2021modeling} on multispecies tumor growth are given.

Different groups prefer to use various finite element method (FEM) libraries, e.g., \texttt{ALBERTA} in \cite{garcke2018multiphase}. \cite{mohammadi2019simulation} utilized element-free Galerkin methods, \cite{wise2008three} a multigrid/finite difference method, and \cite{xu2016mathematical} isogeometric analysis. 
Moreover, the convergence of the FEM in tumor growth has been the subject of theoretical research; see \cite{garcke2022numerical}.

\subsection{Three-dimensional model}
Using the FEM, the 3D models were implemented. The code sequentially solves the system; see Algorithm 2.1 in \cite{fritz2021modeling} for the full model's algorithm. For the potential $\Psi=\Psi_e+\Psi_c$ in the Cahn--Hilliard equation, we employ the classical energy splitting approach, which gives unconditional energy stability; see \cite{elliott1993global}. Thus, the expansive portion $\Psi_e$ is treated explicitly while the contractive portion $\Psi_c$ is treated implicitly. We present the results of numerical experiments in \cite{fritz2021modeling,fritz2021analysis} and demonstrate the relative importance and roles of various biological effects, including cell mobility, proliferation, necrosis, hypoxia, and nutrient concentration, on the generation of MDEs and the degradation of the ECM.

\subsection{Nonlocal phenomena}
Nonlocal effects are not only challenging from an analytical standpoint, but they also pose difficulties for numerical approaches and increase the computational load. The FEM is founded on the notion of local elements, in opposition to the nature of spatial nonlocality. Not only should cells share information within their own element, but also with neighboring elements. In the case of time-fractional PDEs, not only the solution from the previous time step is relevant, but all solutions beginning with the initial condition must also be saved.

\subsubsection{Nonlocal-in-space effects}

In \cite{fritz2019local}, the evolution of the tumor volume fraction was analyzed in both local and nonlocal four-species models.
Thus, we select the gradient-based haptotaxis flux $J_\text{loc}(\phi_V,\ecm)=\chi_h \phi_V \nabla \ecm$ for the local model and $J_\text{nonloc}(\phi_V,\ecm)=\chi_h \phi_V k*\ecm$ for the nonlocal model.
As done in \cite{chaplain2011mathematical}, \cite{gerisch2008mathematical} and \cite{gerisch2010approximation}, we choose a kernel function $k_\eps$, $\eps>0$, in the place of $k$ that approximates the gradient-based haptotaxis effect as $\eps \to 0$.
This also means that a larger nonlocal influence correlates to a greater $\eps$-value.
Specifically, we employ the approximation 
$$\begin{aligned} &(k_\eps*\ecm)(x) - \ecm(x) \cdot (k_\eps*1)(x) \\  &= \int_{\mathbb{R}^d} k_\eps(x-y) (\ecm(y)-\ecm(x)) \dd y
\\&\approx \int_{\mathbb{R}^d} k_\eps(x-y) (\nabla \ecm(x) \cdot (y-x)  ) \dd y
\\&=\nabla \ecm(x) \int_{\mathbb{R}^d} (y-x) \cdot k_\eps(x-y) \dd y  
\\&= \nabla \ecm(x),
\end{aligned}$$
where we selected $k_\eps$ such that $x k_\eps(-x)$ is a Dirac sequence, i.e., it satisfies $\int_{\mathbb{R}^d} x k_\eps(-x) \, \dd x=1$. Specifically, we choose the kernel sequence
$
k_\eps(x)=- \omega(\eps) x \chi_{[0,\eps]} (\vert x \vert_\infty),
$. In the two-dimensional setting, we set the weight $\omega$ depending on $\eps$ to $\omega(\eps)=\frac{3}{8} \eps^{-4}$ in order to fulfill the normalizing Dirac property. %
 .%

\subsubsection{Nonlocal-in-time effects}

We mention the review work by \cite{diethelm2020good} that discusses the pertinent numerical approaches for time-fractional PDEs.
The kernel compressing schemes in \cite{fritz2021equivalence} and \cite{khristenko2021solving}, which reduce the time-fractional PDE to a system of ODEs, are among the numerous efficient methods accessible.
However, the traditional L1 scheme in \cite{oldham1974fractional} is still frequently used due to its simplicity, widespread acceptance, and straightforward implementation, see the survey article by \cite{stynes2021survey}.

Consider the mesh $0=t_0 < t_1 < \dots < t_{N-1}=t_N=T$ of the interval $[0,T]$. The $\alpha$-th Caputo derivative of a given function $\phi$ at $t_n$, $n \in \{1,\dots,N\}$, reads
$$\p_t^{\alpha} \phi(t_{n}) = \frac{1}{\Gamma(1-\alpha)} \int_0^{t_{n}} \frac{\phi'(s)}{(t-s)^\alpha} \ds.$$
We apply the approximation $f'(s)\approx \frac{f(t_{j+1}) - f(t_j)}{t_{j+1}-t_j}$ for $s \in (t_j,t_{j+1})$, which yields
\begin{equation*} \begin{aligned}
    \p_t^{\alpha} \phi(t_{n}) 
    &\approx  \frac{1}{\Gamma(2-\alpha)} \sum_{j=0}^{n-1} w_{n-j-1,n} (\phi(t_{n-j})-\phi(t_{n-j-1})), 
\end{aligned} \end{equation*}
where the weights $w_{m,n}$ for $n,m \in [0,N]$ are given by
$$w_{m,n} = \frac{(t_n-t_m)^{1-\alpha}-(t_n-t_{m+1})^{1-\alpha}}{t_{n-m}-t_{n-m-1}}.$$
The L1 scheme's convergence is in the range of  $\mathcal{O}((\Delta t)^{2-\alpha})$, see \cite{diethelm2020fractional}, and the memory effect is depicted on the right as the history of the preceding time steps $\phi(t_{n-j})$. Exactly this step is computationally intensive due to the need to save the entire history in the computer's memory storage.
Reduce the computational complexity by, for instance, storing only the previous 20 solutions.
Given that the weights on prior solutions drop the further back the previous solution is, this seems reasonable.
But then nothing more can be said about convergence. 

In the works by \citep{fritz2021subdiffusive,fritz2021timefractional} on time-fractional tumor growth models, a fractional linear multistep method is used as in \cite{Lubich86}. Such a method is based on a convolution quadrature scheme, and it generalizes the standard linear multistep method for ODEs. A subclass of these methods generalizes the backward Euler method to fractional settings and approximates the Caputo derivative by
\begin{equation*} \begin{aligned}
\p_t^{\alpha} \phi(t_{n}) 
&\approx  \frac{1}{(\Delta t)^\alpha} \sum_{j=0}^{n-1} (-1)^j \binom{-\alpha}{j} (\phi(t_{n-j})-\phi(0)).
\end{aligned} \end{equation*}
Indeed, setting $\alpha=1$ gives the backward Euler scheme. Similar to the traditional L1 method, it is necessary to store all previous solutions.
The quadrature weights can also be calculated recursively, and such methods are known as Gr\"unwald--Letnikov approximations.  Similar to the classical L1 method, one has to store all the previous solutions, see \cite{diethelm2010analysis} and \cite{Dumitru12} for more details.

\subsection{Uncertainty in tumor modeling}
First off, we note that an orthonormal basis of the Hilbert space $L^2(\Omega)$ on the three-dimensional domain $\Omega=(0,2)^3$ is given by
$$e_{ijk}(x_1,x_2,x_3)=\cos(i\pi x_1/L) \cos(j\pi x_2/L) \cos(k\pi x_3/L),$$
where $L$ is the edge length of the cubic domain $\Omega$.
Then the cylindrical Wiener processes $W_\alpha$ on $L^2(\Omega)$ can be written as
$$W_\alpha(t)(x_1,x_2,x_3) = \sum_{i,j,k=1}^\infty \eta^\alpha_{ijk}(t) e_{ijk}(x_1,x_2,x_3),$$
where $\{\eta^\alpha_{ijk}\}_{i,j,k \in \mathbb{N}}$ is a family of real-valued, independent, and identically distributed (i.i.d.) Brownian motions.
Following the works by \cite{chai2018conforming} and \cite{antonopoulou2021numerical}, we approximate the term involving the Wiener process in the fully discretized system as follows
$$
\frac{1}{\Delta t} \left( \int_{t_n}^{t_{n+1}} \dd W_\alpha(t) , \xi \right)_{L^2(\Omega)} \approx \frac{1}{\Delta t} \sum_{\substack{i,j,k, \\ i+j+k < I_\alpha}} \eta^\alpha_{ijk} (e_{ijk},\xi)_{L^2(\Omega)},
$$
where $\xi \in V_h$ is a test function, $\eta^\alpha_{ijk} \sim \mathcal{N}(0,\Delta t)$ are independent Gaussians, and $I_\alpha$ controls the number of basis functions.

\subsection{Mixed-dimensional coupling}
In the case of 3D-1D tumor growth models, one must implement the new 1D components into the code and establish the link between the 1D and 3D variables.
For time integration of the 1D equations, we employ the implicit Euler method.
For the spatial discretization of the 1D equations, the vascular graph method is used, which corresponds to a node-centered finite volume method, see \cite{reichold2009vascular} and \cite{vidotto2019hybrid} for further details. 

We decouple the 1D and 3D pressure equations at each time step and use block Gau\ss--Seidel iterations to solve the two systems until the 3D pressure converges.
Similarly, the nutrient equation is discretized, with the addition of an upwinding process for the convective term.
The nutrition equations are solved with block Gau\ss--Seidel iterations at each time step.
In \cite{fritz2021modeling}, the numerical approach and discretization of terms that arise in the setting of the 3D-1D coupling are presented in depth.

  \section{Concluding remarks} 
\label{Sec:Conclusion}
We have derived a multiple constituent model from the mass balance law and a Ginzburg--Landau type energy. Like this, we can describe the evolution of tumor cells with various biological phenomena such as angiogenesis. We incorporated stratification and invasion due to ECM deterioration into the model. Moreover, we investigated spatial and temporal nonlocalities, stochasticity resulting from a cylindrical Wiener process, mechanical deformation and elasticity, chemotherapeutic influence, and angiogenesis through mixed-dimensional couplings. 

Like this, we hope that tumor evolution can be studied with all various effects that happen in specific organs. Each tumor is unique, and the parameters have to be tuned for each scenario. One requires a sensitivity analysis with real data and a calibration of the parameters. We regard this as future research after collaborating with doctors and obtaining data.

Mathematically, it cannot be followed immediately whether the nonlinear models are well-posed and admit a solution. There is no unifying theory for the analysis of any nonlinear PDE, and each novel nonlinear system has its own unique challenges that must be examined in depth to confirm or deny the system's well-posedness. We want to emphasize that it is significant to the study the existence of solutions of various models. Otherwise, numerical methods might show solutions, but the model could be ill-posed and not suitable for describing real-world phenomena. 

\backmatter

\bmhead{Acknowledgments}
The author would like to express his sincere thanks to the editor for handling the manuscript.

\bmhead{Data Availability} 
The simulations have been implemented in the code framework “Angiogenesis3D1D” that is accessible on the GitHub project: \url{https://github.com/CancerModeling/Angiogenesis3D1D}.

\section*{Declarations}

\bmhead{Conflict of interest} The author declares no conflict of interest.

\renewcommand*{\bibfont}{\normalfont\scriptsize}
\setlength{\bibsep}{2pt}
\bibliography{literature.bib}
	
\end{document}